\documentclass[preprint,11pt,authoryear]{elsarticle}
\usepackage{amssymb}
\usepackage{amsmath}
\usepackage[margin=1in]{geometry}
\usepackage{enumitem}
\usepackage{url}
\usepackage{float}
\usepackage{subcaption}

\usepackage{amsmath,amssymb,amsfonts}
\usepackage{graphicx}
\usepackage{textcomp}
\usepackage{array}
\usepackage{booktabs}
\usepackage{hyperref}
\usepackage{multirow}
\usepackage{xcolor}

\PassOptionsToPackage{hyphens}{url}\usepackage{hyperref}

\usepackage{tikz}
\usepackage{pgfplots}
\usetikzlibrary{positioning}
\usetikzlibrary{shadows}
\usetikzlibrary{intersections}
\usetikzlibrary{calc}
\usetikzlibrary{spy}
\usetikzlibrary{decorations}
\usepgfplotslibrary{fillbetween}
\usepgfplotslibrary{statistics}
\pgfplotsset{compat=1.3}

\usepackage{algorithm}%
\usepackage{algorithmicx}%
\usepackage{algpseudocode}%
\usepackage{listings}%

\usepackage{comment}
\usepackage{nomencl}
\usepackage{accents}

\algnewcommand\algorithmicinput{\textbf{Input:}}
\algnewcommand\Input{\item[\algorithmicinput]}
\algnewcommand\algorithmicinitialize{\quad\textbf{Initialize:}}
\algnewcommand\Initialize{\item[\algorithmicinitialize]}

\makenomenclature

\renewcommand\nomgroup[1]{%
  \item[\bfseries
  \ifstrequal{#1}{A}{Sets}{%
  \ifstrequal{#1}{B}{Parameters}{%
  \ifstrequal{#1}{C}{Variables}{}}}%
]}

\journal{Applied Energy}

\begin{document}

\begin{frontmatter}



\title{Learning Prosumer Behavior in Energy Communities: Integrating Bilevel Programming and Online Learning } 


\author[DTU]{Bennevis Crowley} 
\author[DTU]{Jalal Kazempour}
\author[DTU]{Lesia Mitridati}
\author[UCSB]{Mahnoosh Alizadeh}

\affiliation[DTU]{organization={Department of Wind and Energy Systems},
            addressline={Elektrovej 325}, 
            city={Lyngby},
            postcode={2800}, 
            country={Denmark}}

\affiliation[UCSB]{organization={Department of Electrical and Computer Engineering, University of California Santa Barbara},
            addressline={Harold Frank Hall}, 
            city={Santa Barbara},
            postcode={93106}, 
            state={CA},
            country={USA}}

\begin{abstract}
Dynamic pricing through bilevel programming is widely used for demand response but often assumes perfect knowledge of prosumer behavior, which is unrealistic in practical applications. This paper presents a novel framework that integrates bilevel programming with online learning, specifically Thompson sampling, to overcome this limitation. The approach dynamically sets optimal prices while simultaneously learning prosumer behaviors through observed responses, eliminating the need for extensive pre-existing datasets. Applied to an energy community providing capacity limitation services to a distribution system operator, the framework allows the community manager to infer individual prosumer characteristics, including usage patterns for photovoltaic systems, electric vehicles, home batteries, and heat pumps. Numerical simulations with 25 prosumers, each represented by 10 potential signatures, demonstrate rapid learning with low regret, with most prosumer characteristics learned within five days and full convergence achieved in 100 days.
\end{abstract}




\end{frontmatter}


\section{Introduction}


The transition to a green energy system requires widespread electrification across various sectors and all levels of society. As many countries and individuals push for this shift, new technologies such as electric vehicles and heat pumps are being increasingly integrated into the energy system. While these technologies, known as Distributed Energy Resources (DERs), help replace fossil fuel-based energy consumption, these novel technologies significantly increase electricity demand. For instance, a single electric vehicle or heat pump can double the annual electricity consumption of a typical residential home \citep{Andersen2017HouseholdConsumption}. Such an increase in electricity consumption comes with challenges. For one, the increased energy costs must be managed to ensure the cost-effective use of these novel assets. Not doing so could lead to catastrophically high costs for some individuals or communities. Such high costs could be avoided through intelligent DER operation, such as load shifting to maximize usage in low-price hours and minimize usage in high-price hours. 

Another challenge arises when considering that these DERs will be connected to the low-voltage distribution grid. As this shift progresses, distribution networks will experience a significant increase in demand, especially during peak-load hours, potentially leading to overloading in local distribution systems. This highlights the need for Distribution System Operators (DSOs) to implement effective congestion management solutions to address these emerging challenges. The inherent flexibility of many DERs can be utilized to alleviate congestion challenges and delay the need for additional distribution grid expansion \citep{Spiliotis2016FlexibilityVsExpansion}. By intelligently adjusting when an electric vehicle charges or a home is heated, these DERs can assist the DSO in operating the distribution system more efficiently and reliably. However, this raises the question of how DERs can be controlled to deliver such demand-side services to the DSO reliably. \cite{Hennig2023CongestionManagementDistribution} provides an overview of possible congestion management services at the distribution level and categorizes them into three main categories: smart network tariffs, market-based approaches, and direct control of flexible loads. Smart network tariffs use price-based signals to manage grid congestion, with static tariffs remaining fixed over the long term and dynamic tariffs adjusting based on system conditions, all aiming to incentivize load shifting to hours when the distribution grid can better handle the demand. Market-based approaches aim to procure flexibility from prosumers—consumers who also produce energy locally—by compensating them for adjusting their load profiles to a desired shape. In contrast, direct control allows the DSO to adjust the consumption of DERs in real-time, addressing congestion issues directly without relying on other actors for congestion management. This paper excludes direct load control due to the high computational burden it places on the DSO \citep{Charbonnier2022} and potential consent issues regarding private consumer assets \citep{Stenner2017}.

To address these issues, this work focuses on utilizing price-based demand response to deliver a service to the DSO by load-shifting the demand of a set of prosumers in an energy community while simultaneously managing the DERs in a cost-effective manner. 

For the sake of generality, we use the term ``prosumers" throughout the paper, although this may also apply to cases where prosumers do not have local production, effectively functioning as simple consumers.

\subsection{Learning prosumer behavior: Why and how?}

To implement price-based demand response consistently, DSOs need access to information on the available flexibility at the distribution level. While this information can be relatively straightforward to obtain for larger assets where demand response is a primary function, it is more complex for many prosumers with DERs. 

The prosumer’s personal usage preferences and habitual behaviors can impose significant restrictions on the flexibility of assets, which a DSO may not know. For instance, the primary function of a heat pump is to maintain the thermal comfort of a building, not to adjust its consumption to support the grid. This creates a challenge for DSOs in accurately predicting the flexibility of such assets. One potential solution is to gather this information directly from prosumers, but frequent data collection could be burdensome, requiring prosumers to repeatedly share details about their habits. Moreover, there is no guarantee that prosumers will be willing to share this information or that they will do so truthfully. An alternative is to \textit{passively learn} this information, thus avoiding the need for continuous or precise communication \citep{Gomez2012LearningPrice}. Passive learning involves observing a prosumer's response to a given signal, such as a price change, to estimate unknown parameters related to behavior. If a large dataset is available for a specific prosumer, supervised learning approaches can be used to model demand response, as demonstrated in \cite{Kim2020} and \cite{Dinh2022}. However, large labeled datasets on individual prosumer  response are unlikely to be readily available. Therefore, methods that do not require extensive labeled data and can quickly learn prosumer behavior through observed responses are preferable in this context. This paper explores the use of online learning techniques to effectively infer prosumer behavior under dynamic pricing.

\subsection{State of the art and research gap}

This literature review is divided into two parts. The first part examines developments in optimal pricing for demand response, while the second part outlines advancements in online learning for demand response.

\subsubsection{Dynamic pricing for demand response}

Dynamic pricing refers to the practice of adjusting prices over time and across different prosumers or locations. Traditionally, residential consumers were subject to fixed prices and tariffs for their electricity usage. However, in recent years, retailers and DSOs have begun offering time-variant prices and tariffs to residential consumers. This shift has shown promising results in various markets. For instance, in 2023 a Danish DSO achieved a 10\% load shift by implementing higher tariffs in peak-load hours \citep{JanuarFebruar2023}. Meanwhile, as early as 2010, \cite{Faruqui2010} demonstrated a demand response of up to 44\% from residential users when provided with the right enabling technologies. These figures are expected to increase with the growing adoption of automated home energy management systems and enhanced flexibility in the distribution grid.

Despite these advances, the key challenge remains determining the optimal way to set dynamic prices that effectively stimulate the desired demand response. Some studies have explored pricing electricity based on the operations of the distribution grid, using schemes like distribution locational marginal pricing \citep{Sotkiewicz2006DLMP,Papavasiliou2018,Mieth2020}, where prices vary both temporally and spatially. These approaches have been expanded to include price-based demand response, such as in \cite{PANDEY2025111398}, which demonstrate that differentiated pricing can be an effective tool in managing DERs.


A common approach for setting dynamic prices while considering prosumer responses is the Stackelberg (leader-follower) game. In this framework, a price setter, acting as the leader of the game, endogenously accounts for the follower's optimization problem while optimizing their own objective function. The use of the Stackelberg game in electricity pricing problems dates back to \cite{Luh1982}, which first examined the interaction between retailers and consumers responding to prices. Since then, the application of Stackelberg games has expanded significantly. For example, \cite{jia2016dynamic} provides an in-depth economic analysis of a retailer setting day-ahead dynamic prices for consumers through a Stackelberg game, investigating the welfare impacts of integrating renewable energy sources. \cite{Alves2020} focuses on optimizing both the prices and the periods of retailer electricity prices, considering consumer cost minimization. \cite{Venkatraman2022} applies the same method to set optimal distribution grid tariffs for a DSO, evaluating the effects on electric vehicle charging. \cite{Wei2022} extends the scope of the follower's problem to include demand response from multi-energy buildings, while the retailer continues to set prices as the leader of the game. \cite{Luo2023} further explores multi-energy systems, incorporating energy storage systems as followers as well. Lastly, \cite{Taherian2021} integrates consumer behavior predictions into a bilevel model, while still focusing on optimal price setting from the retailer's profit-maximizing perspective.

Although various techniques have been applied to solve such bilevel problems, all of these studies rely on reformulation or decomposition methods to adopt a single-shot optimization approach. These models do \textit{not} learn about prosumer behavior and flexibility over time. Instead, they either assume the availability of a large dataset to train a predictive model or adopt a simplified price-response model with \textit{known} parameters of the follower, which may not accurately reflect actual prosumer behavior. Failing to update this assumed knowledge can be problematic, as inaccurate predictions or incorrect assumptions can lead to suboptimal results in the bilevel model.

\subsubsection{Online learning for demand response} \label{subsubsec:literature online}

In recent years, numerous studies have explored the application of online decision-making techniques to demand response, many of them focusing on the bandit framework. The bandit framework is a model for sequential decision-making, where an agent must choose an action when faced with uncertain rewards or costs. After each action, the agent receives a reward but does not know the reward associated with other possible actions. This sequence is repeated over several rounds, where a new decision is made every round based on the information gathered up until that point. The goal within this framework is to maximize the cumulative reward over time, which requires balancing exploration (trying out new actions) and exploitation (using known actions that yield higher rewards). This tension between exploration and exploitation is known as the \textit{exploration-exploitation trade-off}.


Several algorithms are available to navigate this trade-off, including greedy, $\epsilon$-greedy, Upper Confidence Bound (UCB), and Thompson sampling. While these algorithms differ in how they select actions, they all aim to maximize the agent's reward over a series of decisions. The simplest approach, the greedy algorithm, always chooses the action with the highest expected reward without considering exploration. The $\epsilon$-greedy algorithm extends this by selecting a random action with a small probability in each round, thus enabling exploration. The UCB algorithm also aims to select the action with the highest expected reward but it adds an optimistic boost to the expected reward based on the uncertainty of an action’s potential reward (i.e., how little it has been explored). This boost encourages the exploration of less-known actions. Of course, many other algorithms can be used in a bandit setting.

The bandit framework has been regularly applied in the demand response setting for a bit over a decade at this point. For example, \cite{Taylor2014Index} formulates the peak-shaving demand response problem within the bandit framework, aiming to select the optimal set of consumers for load curtailment. This work derives index policies for the demand response formulation and demonstrates that these policies outperform a greedy algorithm. Similarly, \cite{Rajagopal2015OnlineLearning}  seeks to select the optimal set of demand response customers but employs a different bandit formulation, for which an optimal learning algorithm is provided. In a related study, \cite{Li2020} applies a combinatorial version of the UCB  algorithm to learn the optimal subset of consumers needed to achieve a specific amount of demand response. 

Moving beyond consumer selection and opt-outs, \cite{Schneider2022} uses a UCB-based algorithm to estimate consumer baseline profiles over time while simultaneously learning to offer more effective incentive payments for demand response. Bandit frameworks have also been applied to model the demand response of thermostatically controlled loads. For example, \cite{Antoine2017TCL} utilizes an adversarial bandit framework to learn the parameters of such loads based on their response to direct load control signals. This approach is extended in \cite{Amr2017TCLs}, where the bandit framework is used to optimize the response of thermostatically controlled loads, enabling the selection of the best candidates for ancillary service provision to the grid.

Looking at bandit applications for price-based demand response, \cite{jia2014onlinelearningapproachdynamic} uses a bandit framework to adjust dynamic prices, learning the parameters of a linear price response function for consumers. Similarly, \cite{Khezeli2018} examines an electric power utility that adjusts its prices in a risk-sensitive manner to reduce demand. The study employs continuous least-squares estimation to learn the assumed affine price sensitivity functions of a set of consumers. \cite{Mieth2020Online} extends this work by using least squares estimation to learn a quadratic consumer price sensitivity function while incorporating network constraints through distributionally robust power flow. Finally, \cite{Hutchinson2024} explores the safe operation of safety-critical networks, such as power grids, while maximizing consumer utility. This study uses an algorithm based on UCB to learn consumers' price response over time.

This paper focuses on the last of the aforementioned algorithms, Thompson sampling, chosen for its superior performance compared to greedy algorithms and its ability to handle more complex problems than UCB-based methods. Thompson sampling adopts a Bayesian approach to multi-armed bandits, where an initial distribution represents the agent's knowledge of the unknown parameters. This distribution is updated across all rounds to reflect new information. In each round, a sample is drawn from this distribution and assumed to represent the true value of the unknown parameter. By sampling, the exploration-exploitation trade-off is inherently balanced by the algorithm. 

Thompson sampling has been successfully applied to demand response problems in other research. For example, \cite{Chen2021Learning} and \cite{Chen2021Contextual}   use Thompson sampling to examine how user preferences and environmental factors influence consumer opt-out decisions in a direct load control demand response scenario. \cite{Moradipari2018Price} applies Thompson sampling to set real-time prices for an aggregator seeking to control the loads within its portfolio while learning the underlying price-response structure. \cite{Tucker2020ConTS} extends this approach by aligning an aggregated set of demands with a daily load profile while simultaneously learning an unknown parameter vector that characterizes the price-based load response of the aggregated demands, ensuring safe operation under the given price signals.

While many of these studies optimize price choices within the algorithm, none of the identified pricing approaches incorporate optimization problems that reflect the leader-follower structure of price-based demand response in the price-setting process. Incorporating this leader-follower structure to create a more refined price-setting strategy could lead to more accurate optimal price signals for prosumers, hereby enhancing demand response outcomes. Additionally, the pricing strategies used in these studies typically target clusters of appliances rather than individual consumers. Implementing individual-level pricing could allow for more effective demand response from price-sensitive prosumers.


\subsection{Contributions}

Given the above literature review and the research gaps identified, the main contributions of this paper include:

\begin{itemize} 
\item Bilevel optimization is integrated with an online learning framework using Thompson sampling, enabling the price setter to learn and adapt to prosumer behavior over time. To the best of our knowledge, this is the first paper to apply this combination, offering a novel and flexible approach to dynamic price-setting in demand response. 

\item Unlike traditional single-shot bilevel models, this work removes the assumption of perfect knowledge about prosumer behavior. Instead of relying on historical data-driven prediction models, the price setter learns and adapts to prosumer responses in real-time. 

\item In contrast to prior works that use Thompson sampling for aggregate-level learning, this paper focuses on individual-level learning, allowing for more tailored and effective pricing strategies.

\end{itemize}

Without loss of generality, we implement our proposed method within the context of an \textit{energy community} consisting of a set of prosumers. The price setter (leader) in this case is the community manager, a non-profit entity, who employs our integrated algorithm to set dynamic prices for prosumers (followers) over a specified time period, while simultaneously learning their behavior by observing their responses. The community manager’s objective is to provide grid services to the DSO by capping the community’s total consumption during peak-load hours, thereby enabling the community to benefit from discounted grid tariffs. We evaluate the performance of the proposed method using a regret metric. The proposed method can also be applied to any demand-response application in which an aggregation of prosumers, coordinated through dynamic pricing, collectively offers any form of flexibility services.

\subsection{Paper organization}

The remainder of this paper is organized as follows: Section \ref{sec:framework} presents the proposed framework, detailing how we integrate online learning with bilevel optimization. Section \ref{sec:math} provides an in-depth discussion of the mathematical formulation of the framework, the implemented algorithm, and the performance metrics used to evaluate the algorithm. Section \ref{sec:case} presents a numerical case study to analyze the performance of the framework in a realistic setting. Finally, Section \ref{sec:conclusion} summarizes the key findings, explores the potential impact on distribution-level congestion management, and outlines directions for future research in this area.

\textit{Notation}: A full nomenclature can be found in Appendix A, but the general notation is introduced here. Sets are denoted by upper-case calligraphic letters, e.g., $\mathcal{T}$ and $\mathcal{I}$. Parameters are represented by upper-case letters, such as $\mathrm{L}, \mathrm{PV}$, or Greek letters, like $\alpha$, and may be indexed by sets when necessary. Single variables are written in lowercase italics with the relevant indices, e.g., $x_{n,t}$, while vectors of variables are written in bold lowercase with an index if required, e.g., $\boldsymbol{x}_t$. Matrices of variables are written in bold upper-case letters, e.g., $\boldsymbol{X}$. Some sets, parameters, and variables may have superscripts to distinguish them from similar terms. Parameters and variables that are sampled from a distribution, or result from an optimization based on a sampled value, are denoted with a tilde, i.e., $\tilde{(\cdot)}$. Any deviations from this notation will be clarified in the text.
\section{Proposed framework: Integrating learning and optimization} \label{sec:framework}

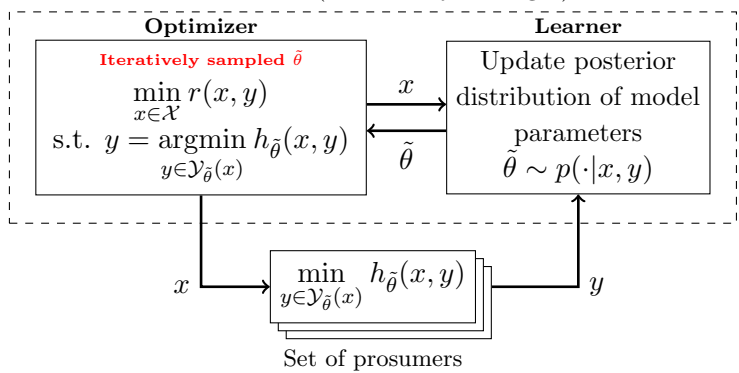
\begin{figure}[ht]
    \centering
    \begin{subfigure}[b]{0.33\textwidth}
        \centering
        \begin{tikzpicture}[align = center]
    \node[rectangle, draw = black, minimum width = 125pt, minimum height = 60pt] (Math) {\vspace{30pt} \:{\color{red}{\tiny{\textbf{Assumption: Perfect knowledge of $\theta$}}}}\: \\ $\underset{x \in \mathcal{X}}{\min} \: r(x,y)$ \\ s.t. $y = \underset{y \in \mathcal{Y}_{\theta}(x)}{\text{argmin}} \: h_{\theta}(x,y)$ \\};
    \node[above = 0pt of Math.north] (box text) {\small{Price setter (community manager)}};
    \node[rectangle, draw = black, below = 30pt of Math.south, fill = white, copy shadow = {shadow xshift = 6pt, shadow yshift = -6pt}, copy shadow = {shadow xshift = 3pt, shadow yshift = -3pt}] (Environment) {$\underset{y \in \mathcal{Y}_{\theta}(x)}{\min} \: h_{\theta}(x,y)$};
    \node[below = 6pt of Environment.south] (environment text) {\footnotesize{Set of prosumers}};

    \draw[->, line width = 1pt] ([shift = {(-5pt,0pt)}]Math.south) -- ([shift = {(-5pt,0pt)}]Environment.north) node[midway, left = 0pt] {$x$};
    \draw[->, line width = 1pt] ([shift = {(5pt,0pt)}]Environment.north) -- ([shift = {(5pt,0pt)}]Math.south) node[midway, right = 0pt] {$y$};
\end{tikzpicture}
        \caption{\scriptsize{General bilevel programming \\ (available in the literature)}}
        \label{fig:bilevel}
    \end{subfigure}
    \begin{subfigure}[b]{0.63\textwidth}
        \centering
        \begin{tikzpicture}[align = center]
    \node[rectangle, dashed, draw = black, minimum width = 275pt, minimum height = 80pt] (Agent) {};
    \node[above = 0pt of Agent.north] (agent text) {\small{Price setter (community manager)}};
    \node[rectangle, draw = black, minimum width = 125pt, right = 10pt of Agent.west] (Optimizer) {{\color{red}{\tiny{\textbf{Iteratively sampled $\tilde\theta$}}}}\\ $\underset{x \in \mathcal{X}}{\min} \: r(x,y)$ \\ s.t. $y = \underset{y \in \mathcal{Y}_{\tilde\theta}(x)}{\text{argmin}} \: h_{\tilde\theta}(x,y)$};
    \node[above = -3pt of Optimizer.north] (optimizer label) {\scriptsize{\textbf{Optimizer}}};
    \node[rectangle, draw = black, minimum width = 2cm, minimum height = 1cm, left = 10pt of Agent.east] (Learner) {\small{Update posterior} \\ \small{distribution of model} \\ \small{parameters} \\ $\tilde{\theta} \sim p (\cdot \vert x,y) $};
    \node[above = -1pt of Learner.north] (learner label) {\scriptsize{\textbf{Learner}}};
    \node[rectangle, draw = black, below = 10pt of Agent.south, fill = white, copy shadow = {shadow xshift = 6pt, shadow yshift = -6pt}, copy shadow = {shadow xshift = 3pt, shadow yshift = -3pt}] (Environment) {$\underset{y \in \mathcal{Y}_{\tilde\theta}(x)}{\min} \: h_{\tilde\theta}(x,y)$};
    \node[below = 6pt of Environment.south] (environment text) {\footnotesize{Set of prosumers}};

    \draw[->, line width = 1pt, draw = black] (Optimizer.south) |- (Environment.west) node[midway, left = 0pt] {$x$};
    \draw[->, line width = 1pt, draw = black] ([shift={(6pt,0pt)}]Environment.east) -| (Learner.south) node[midway, right = 0pt] {$y$};
    \draw[->, line width = 1pt, draw = black] ([shift = {(0pt,-5pt)}]Learner.west) -- ([shift = {(0pt,-5pt)}]Optimizer.east) node[midway, below=0pt] {$\tilde{\theta}$};
    \draw[->, line width = 1pt, draw = black] ([shift = {(0pt,5pt)}]Optimizer.east) -- ([shift = {(0pt,5pt)}]Learner.west) node[midway, above=0pt] {$x$};
    \end{tikzpicture}
        \caption{\scriptsize{Proposed integrated framework:  A combination of bilevel programming and online learning (Thompson sampling)}}
        \label{fig:contribution}
    \end{subfigure}
    \caption{Comparison of state-of-the-art bilevel framework and the proposed integration of the bilevel problem into an online learning framework}
    \label{fig:comparison}
\end{figure}

To further clarify our contributions in a more mathematical sense, while still maintaining the general perspective, we compare the general bilevel optimization problem shown in Figure \ref{fig:bilevel} with our proposed integrated method (a combination of bilevel programming and online learning), illustrated in Figure \ref{fig:contribution}. In both figures, the leader is the price setter (community manager), who sets dynamic prices $x \in \mathcal{X}$ in the upper level of the bilevel structure. The followers in the lower level are the prosumers (community members), who respond to prices $x$ by adjusting their consumption levels $y \in \mathcal{Y}(x)$. In both cases, the price setter optimizes its objective function $r(x,y)$, which depends on both prices $x$ and responses $y$. For given prices $x$, each prosumer in both figures minimizes its own objective $h(x,y)$ by optimizing $y$. In both figures, as indicated by the bilevel structure, the price-setter's optimization problem is constrained by the optimization problem of each prosumer.

What differentiates these two figures is whether the leader has perfect knowledge of the prosumers' assets and how they prefer to use them, which we refer to as \textit{signatures}. In Figure \ref{fig:comparison}, this information is represented by $\theta$. These signatures may represent different types of prosumers, such as those who are fully or partially flexible in shifting load throughout the day or those who possess specific assets like Photovoltaic (PV) systems, Electric Vehicles (EVs), heat pumps, etc. For example, one type of prosumer may be flexible only during certain hours of the day and possess a PV system but not an EV or a heat pump.  In Figure \ref{fig:bilevel}, it is assumed that the community manager has perfect knowledge of the prosumers' assets and how the prosumers use them (e.g., what their EV charging needs are), which is a strong assumption. In contrast, in the proposed framework, illustrated in Figure \ref{fig:contribution}, the community manager employs an online learning method to estimate this information over time, represented by $\tilde{\theta}$, by observing their responses $y$ to the prices $x$. Acting as the learner, the community manager updates the posterior distribution of the model parameters $\tilde{\theta} \sim f (\cdot \vert x,y)$ after each observation and draws a sample for the next decision-making point. Thus, our model does not rely on perfect knowledge of which prosumer adheres to which signatures; instead, it learns this information iteratively while optimally setting dynamic prices based on the available information at each decision point.

We begin by providing an overall description of the energy community, its grid service to the DSO, and the objectives and interactions of both the community manager and the prosumers. Next, we explain how the signatures of prosumers are incorporated into the bilevel model. Finally, we present the general framework of our proposal for integrating bilevel programming and online learning, which will be applied to dynamic pricing in the energy community.


\subsection{Energy community: Grid service, objectives, and interactions}

The goal of the price-setting problem for the community manager is to shape demand response so that the total energy consumption at the grid connection point remains below a specified threshold. This threshold is not necessarily a physical capacity limit, but rather a level agreed upon with the DSO, which can be considered a grid service provided by the community, referred to as \textit{capacity limitation services}. This service ensures that the DSO can manage the load from the energy community, preventing congestion in the upstream grid that could affect other connected agents. In return for adhering to the capacity limitation, the community benefits from discounted distribution grid tariffs. However, exceeding the threshold results in a penalty from the DSO to cover associated costs.

Thus, the community manager’s objective is to minimize the total electricity costs for the community, considering daily power imports $\mathbf{p^{\mathrm{im}}}$, power exports $\mathbf{p^{\mathrm{ex}}}$, and any penalties incurred for exceeding the capacity limitation. This is achieved by setting time- and prosumer-specific prices $\mathbf{X}$ that influence prosumers' consumption. In addition to cost minimization, the community manager must ensure the feasibility of the resulting power flows. This includes satisfying power balance constraints at each node, enforcing bounds on power imports/exports and prices, and meeting economic constraints, such as budget balance for the community manager and individual rationality for each prosumer. 

Each prosumer seeks to minimize electricity procurement costs based on the prices set by the community manager, while satisfying constraints related to consumption needs, flexibility, and preferences. Ideally, these individual problems could be embedded within the community manager’s price-setting model. However, this requires accurate and truthful communication of all relevant input parameters, which is often impractical in real-world scenarios.

To address this challenge within the bilevel optimization framework, the community manager can use learning mechanisms to infer prosumers' constraints and preferences based on their observed responses to prices, along with prior knowledge of their assets, all without interfering with the community’s regular operation. Consequently, the formulation of the lower-level problems must be adjusted, as detailed in the following subsection.

\subsection{How to incorporate signatures into the bilevel model?}

To continue our high-level discussion of the model used in this paper, we now move to the prosumers demand model. Demand profiles at the residential level are typically observed at an aggregate household level, combining the power usage of multiple different assets within the household. Modeling the aggregate profile without considering the individual assets may overlook critical asset-level constraints that are influenced by prosumer preferences or behavior. Therefore, to accurately model a prosumer's demand profile, it is essential to disaggregate the demand data to the asset level.

A straightforward approach to disaggregation is to create a mathematical model for each asset, assuming that the exact constraints and prosumer preferences are known. However, there may be multiple usage scenarios for any given asset. For instance, an EV might be disconnected at different times, depending on when the prosumer needs to drive, which in turn affects the EV's availability to provide grid services. These various usage scenarios cannot be fully captured in a single model. To address this, we develop separate mathematical models for each potential usage case, such as an EV being disconnected early in the morning or in the evening. Each of these mathematical models represents a  \textit{signature}. For each prosumer $n$, there are multiple signatures $k$ where each signature corresponds to a distinct daily demand or production profile, denoted as $\mathbf{p}_{n,k}$. For a given prosumer, the collection of all these profiles is defined as 
\begin{equation}
    \mathbf{P}_n = \left[ \mathbf{p}_{n,k} \right]_{\forall k \in \mathcal{K}}.
\end{equation}

The observed aggregate daily response $\boldsymbol{y}_n$ for prosumer $n$ can be obtained by weighting each signature $k$ by a parameter $\theta_{n,k}$ and summing the resulting contributions. For instance, consider an EV: if the prosumer disconnects their car in the morning, the signature corresponding to the morning disconnection will have a non-zero weight, while the evening signature will have a zero weight, as it does not appear in the prosumer's aggregate demand profile. The signatures contributing to the aggregate response will vary across prosumers, which is reflected in the differing weights for each signature across all prosumers. Furthermore, demand curves exhibit some inherent variability on a daily basis, which can be modeled by introducing a noise term $\boldsymbol{\epsilon}$. This leads to the linear prosumer response model, as shown in \eqref{eq:response_model}, where the signature weights of each prosumer $\boldsymbol{\theta}_n$ are unknown parameters to be learned:

\begin{equation} \label{eq:response_model}
    \boldsymbol{y}_n = \mathbf{P}_n \boldsymbol{\theta}_n + \boldsymbol{\epsilon}.
\end{equation}

In the context of this work, the energy community manager aims to learn the unknown weights of the signatures for all prosumers and their corresponding signatures. Each day, the community manager optimizes dynamic prices based on the best available estimate of the community's unknown weights. In response, the prosumers adjust their demand according to their true underlying weights, generating a load response that the community manager observes. A detailed mathematical discussion of the signatures is provided in Section \ref{subsubsec:consumer}

\begin{figure}[t]
    \centering
    \input{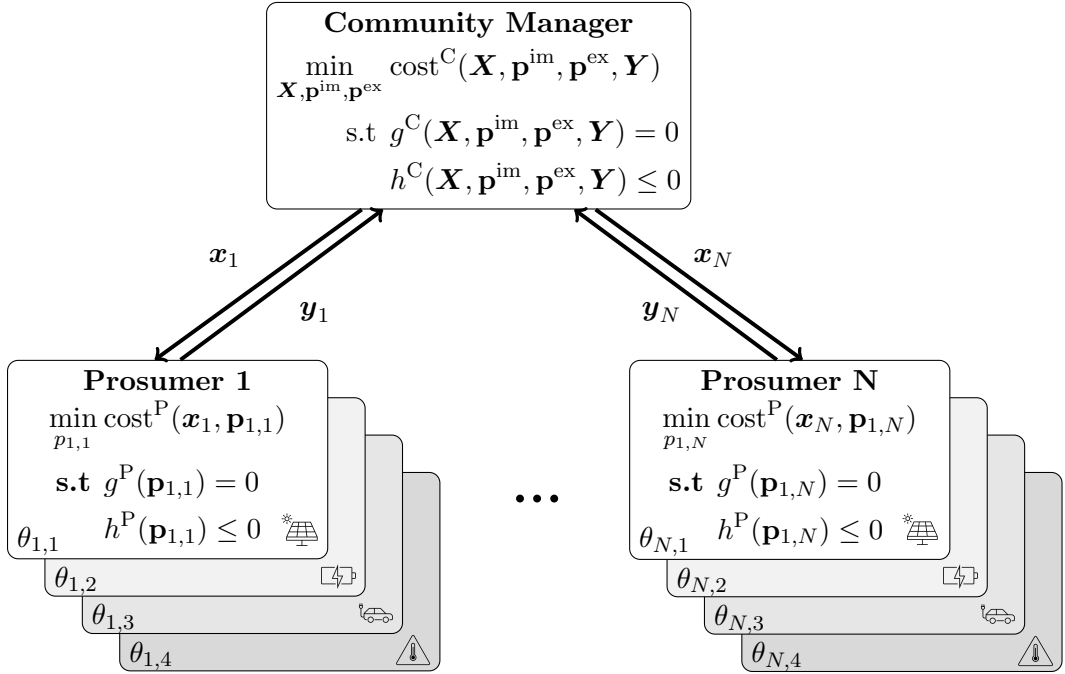}
    \vspace{3ex}
    \caption{Illustration of the bilevel problem setting dynamic prices in the energy community: Four signatures have been arbitrarily selected for each prosumer for illustration purposes. These four signatures correspond to four different types of assets: PV, battery, EV, and heat pump. The number of signatures is not fixed and can be adjusted based on the needs of the community manager.}
    \label{fig:new-bilevel}
\end{figure}

Figure \ref{fig:new-bilevel} illustrates the structure of the proposed bilevel model, incorporating signatures $k = \{1, 2, 3, 4\}$ for each prosumer $n = \{1, ..., N\}$. We have arbitrarily selected four signatures for illustration purposes, but more signatures will be considered in our case study. While the community manager in the upper level cannot directly observe the true weights $\boldsymbol{\theta}^*_{n}$ used to generate responses $\boldsymbol{y}_n$, it does have access to the metered response data. This observable data can then be used to update the manager's estimate of the weights $\hat{\boldsymbol{\theta}}_n$. Since a prosumer's dynamic prices $\boldsymbol{x}_n$ must be shared on a daily basis, the community manager cannot wait for a large batch of data to improve its estimates. As a result, the manager must learn in an online manner, updating the estimates as new data becomes available.

\subsection{Thompson sampling: Learning prosumer behavior by updating $\theta_{n,k}$}




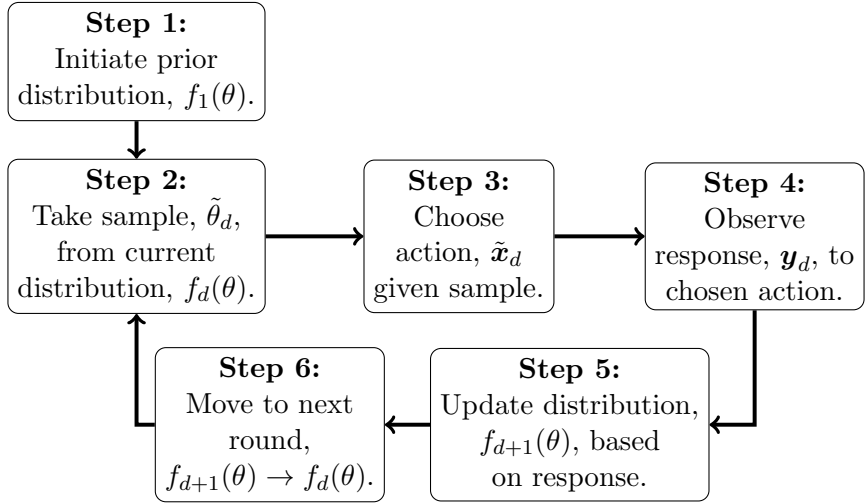
\begin{figure}[t]
    \centering
    \begin{tikzpicture}[align = center, 
                    arrow/.style={->,line width = 1.5pt}, 
                    box/.style={draw, rectangle, rounded corners, minimum width=2.5cm, minimum height=1cm, align=center},
                    ellipse/.style={draw, ellipse, minimum width=2.5cm, minimum height=1cm, align=center}]
        
    \node[box] (prior) {\textbf{Step 1:} \\ Initiate prior \\distribution, $f_1(\theta)$.} ;
    \node[box, below = 0.5cm of prior] (sampling) {\textbf{Step 2:} \\ Take sample, $\tilde{\theta}_d$,\\ from current \\distribution, $f_d(\theta)$.}; 
    \node[box, right = 1.3cm of sampling](action) {\textbf{Step 3:} \\ Choose \\ action, $\tilde{\boldsymbol{x}}_d$ \\ given sample.}; 
    \node[box, right = 1.2cm of action] (response) {\textbf{Step 4:} \\ Observe \\ response, $\boldsymbol{y}_d$,  to \\ chosen action.}; 
    \node[box, below left = 0.5cm and 0.6cm of response.south] (update) {\textbf{Step 5:} \\ Update distribution, \\ $f_{d+1}(\theta)$, based \\ on response.} [below = 2cm of action];
    \node[box, left = 0.6cm of update] (next) {\textbf{Step 6:} \\ Move to next \\round, \\ $f_{d+1}(\theta) \rightarrow f_{d}(\theta)$.};
    
    \draw[arrow] (prior) -- (sampling);
    \draw[arrow] (sampling) -- (action);
    \draw[arrow] (action) -- (response);
    \draw[arrow] (response.south) |- (update.east);
    \draw[arrow] (update.west) -- (next.east);
    \draw[arrow] (next.west) -| (sampling.south);

\end{tikzpicture}
    \caption{Illustration of the iterative procedure for Thompson sampling}
    \label{fig:ts-framework}
\end{figure}

We use Thompson sampling as an online learning algorithm \citep{russo2020tutorialthompsonsampling}. Below, we provide an intuitive introduction to the Thompson sampling algorithm, focusing on a single unknown parameter. The algorithm begins by setting a prior distribution $f_1(\theta)$ for the unknown parameter of interest, $\theta$. It updates this distribution throughout the learning process to effectively estimate the true underlying parameter $\theta^*$ over time. In each iteration $d\!\in\!\mathcal{D}$,  representing a day in our case study, a sample $\tilde\theta_d$ is drawn from the current distribution $f_d(\theta)$. This sample is assumed to represent the true value of the unknown parameter. Dynamic prices on that day $\tilde{\boldsymbol{x}}_d$ are determined based on this sample. The community manager then observes the true response of prosumers $\boldsymbol{y}_d$ to these dynamic prices. These observations are used to update the prior distribution via Bayes' rule, resulting in a posterior distribution. For the next day, this posterior distribution becomes the new prior from which the sample is drawn. This iterative procedure of sampling from and updating the posterior distribution is summarized in Figure \ref{fig:ts-framework}.

One appealing characteristic of Thompson sampling is that it inherently balances exploration and exploitation due to its Bayesian nature. As the algorithm samples from the posterior distribution, new parameter values are naturally generated and tested over the iterative rounds, facilitating exploration. Over time, as sufficient information is gathered from the posterior distribution, the algorithm will sample parameters closer to their true values, yielding higher rewards and enabling exploitation. Therefore, no additional metrics or hyperparameters are required to allow the algorithm to balance exploration and exploitation.

Furthermore, any prior distribution can be used to model the unknown parameter in Thompson sampling. However, the choice of prior distribution and the relationship between the observed reward and the unknown parameter can significantly affect the complexity of the posterior update. This paper focuses on combining the Thompson sampling algorithm with a more sophisticated optimization model to select the optimal actions at each iteration, rather than developing new insights into the algorithm itself. Thus, we intentionally choose simple yet appropriate prior distributions and reward functions to model the response of prosumers.

To assess the performance of Thompson sampling, we focus on \textit{regret}. Regret is defined as the cumulative difference in reward between the algorithm in use and an optimal algorithm with full knowledge of the underlying parameters. This comparison quantifies the reward lost due to the lack of complete information. A well-performing learning algorithm aims to minimize regret over time, ideally achieving near-zero regret in later iterations, which indicates that the true underlying parameters of the model have been successfully learned.


The proposed integration of optimization and learning approaches in this paper combines the bilevel problem for determining dynamic prices with the Thompson sampling algorithm. To clarify this integration, each step in Figure \ref{fig:ts-framework} corresponds to a specific action:

\begin{itemize}

\item \textbf{Step 1:} Initialize a prior distribution for the signature weights of all prosumers. 

\item \textbf{Step 2:} On a given day, a sample is drawn from each distribution and assumed to represent the true signature weight for that day. 

\item \textbf{Step 3:} 
Use the bilevel problem to determine the optimal prices that maintain the capacity limitation and minimize the total community cost. The samples drawn from each distribution are incorporated into the price-setting process. 

\item \textbf{Step 4:} Prosumers in the community adjust their consumption according to their true signature weights. The  community manager observes the response of each prosumer via metering data. 

\item \textbf{Step 5:} The  community manager updates each distribution based on the observed responses from all prosumers. 

\item \textbf{Step 6:} The updated distribution is used as the prior for the next day. The algorithm then returns to \textbf{Step 2}. 
\end{itemize}

The detailed mathematical formulation of both the bilevel model and the proposed Thompson sampling algorithm, along with its evaluation metrics, are presented in the following section.


\section{Mathematical formulation} \label{sec:math}

This section delves into the mathematical details of the bilevel model with signatures, as illustrated in Figure \ref{fig:new-bilevel}, followed by the Thompson sampling algorithm, depicted in Figure \ref{fig:ts-framework}. We eventually introduce the regret metric used to evaluate the algorithm's performance.

\subsection{Bilevel optimization problem with signatures}

We begin by presenting the upper-level problem of the community manager, followed by the lower-level problems that represent the prosumers' signatures.

\subsubsection{The upper-level problem: Dynamic price setting}

The upper-level problem aims to minimize the total community cost by determining optimal prices for prosumers across hours, while also providing capacity limitation services to the DSO:

\begingroup
\begin{subequations} \label{eqs:upper}
\begin{alignat}{2}
    \underset{\Xi^{\rm{U}}}{\operatorname{min}} \ &\sum_{t \in \mathcal{T}} \Bigg[ p^{\mathrm{im}}_t \Big(\rho^{\mathrm{spot}}_t + \gamma^{\mathrm{im}}_t\Big) - p^{\mathrm{ex}}_t \Big(\rho^{\mathrm{spot}}_t &&- \gamma^{\mathrm{ex}}_t\Big) + \alpha^{\mathrm{DSO}}_t p^{\mathrm{pen}}_t \Bigg] \label{eq:upper_obj} \\ 
    \hspace{0.25cm} \mathrm{s.t.} \: &p^{\mathrm{im}}_t - p^{\mathrm{ex}}_t = \sum_{n \in \mathcal{N}} y_{n,t}  && \forall t \in \mathcal{T} \label{eq:upper_balance}\\
    \hspace{0.25cm} \hspace{0.59cm} &p^{\mathrm{pen}}_t \geq p^{\mathrm{im}}_t - \overline{\mathrm{P}}^{\mathrm{DSO}}_t && \forall t \in \mathcal{T} \label{eq:upper_pen}\\
    \hspace{0.25cm} \hspace{0.59cm} &p^{\mathrm{im}}_t, p^{\mathrm{ex}}_t, p^{\mathrm{pen}}_t \geq 0 && \forall t \in \mathcal{T} \label{eq:upper_pos}\\
    \hspace{0.25cm} \hspace{0.59cm} &y_{n,t} = \sum_{k \in \mathcal{K}} p_{n,t,k} \theta_{n,k} \  &&\forall n \in \mathcal{N}, t \in \mathcal{T} \label{eq:upper_response} \\
    \hspace{0.25cm} \hspace{0.59cm} &\sum_{t \in \mathcal{T}}x_{n,t} y_{n,t} \leq \text{cost}_n^{\mathrm{ex}} &&\forall n \in \mathcal{N} \label{eq:upper_IR}\\
    \hspace{0.25cm} \hspace{0.59cm} &\sum_{n \in \mathcal{N}} \sum_{t \in \mathcal{T}}x_{n,t} y_{n,t} \geq \sum_{t \in \mathcal{T}}p^{\mathrm{im}}_t \Big(\rho^{\mathrm{spot}}_t +&& \gamma^{\mathrm{im}}_t\Big) - p^{\mathrm{ex}}_t \Big(\rho^{\mathrm{spot}}_t - \gamma^{\mathrm{ex}}_t \Big) + \alpha^{\mathrm{DSO}}_t p^{\mathrm{pen}}_t, \label{eq:upper_RA}
\end{alignat}
\end{subequations}
\endgroup
where $\Xi^{\rm{U}} = \{ p^{\mathrm{im}}_t , \ p^{\mathrm{ex}}_t, \ p^{\mathrm{pen}}_t, \ x_{n,t}, \ y_{n,t}\}$ is the set of upper-level decision variables. The upper-level objective function \eqref{eq:upper_obj} minimizes the total cost of the energy community. The first term represents the cost of importing power, where the imported power $p^{\mathrm{im}}_t$ at the connecting point to the distribution grid at time $t \in \mathcal{T}$ is multiplied by the day-ahead price $\rho^{\mathrm{spot}}_t$ plus any grid tariffs $\gamma^{\mathrm{im}}_t$. The second term in \eqref{eq:upper_obj} reflects the revenue from exporting power to the grid, where the exported power $p^{\mathrm{ex}}_t$ is multiplied by the day-ahead price, minus any export tariffs $\gamma^{\mathrm{ex}}_t$. The final term accounts for the cost of violating the capacity limitation imposed by the DSO. This term multiplies the excess power $p^{\mathrm{pen}}_t$ by the pre-defined penalty $\alpha^{\mathrm{DSO}}_t$ to be paid to the DSO. To prioritize adherence to the capacity limitation, the penalty for exceeding the limit is set significantly higher than the costs or revenues associated with importing or exporting power.

The objective function \eqref{eq:upper_obj} is subject to physical constraints governing the operation of the grid within the community. Constraint \eqref{eq:upper_balance} ensures that the total power imported and exported is balanced with the aggregate power consumption and production of all prosumers $n \in \mathcal{N}$.  Constraint \eqref{eq:upper_pen} tracks the amount of imported power that exceeds the pre-defined capacity limit $\overline{P}^{\mathrm{DSO}}_t$, ensuring that any violations are accounted for at all times. Constraint \eqref{eq:upper_pos} enforces non-negativity of the variables. Constraint \eqref{eq:upper_response} calculates the expected response of prosumers by multiplying the signature weights $\theta_{n,k}$ by the power signatures $p_{n,t,k}$  from each individual prosumer. This is equivalent to \eqref{eq:response_model}, though it is not shown here in matrix form. Finally,  \eqref{eq:upper_IR} and \eqref{eq:upper_RA} ensure individual rationality for each prosumer and revenue adequacy for the community manager, respectively. Individual rationality is maintained by ensuring that the cost paid by each prosumer $n$ for its total power imports and exports is less than the cost the prosumer would incur if operating independently outside the community, $\text{cost}_n^{\mathrm{ex}}$. Revenue adequacy ensures that the community manager collects sufficient payments from prosumers to cover any obligations to external stakeholders (e.g., the DSO, retailer, etc.), thereby guaranteeing that the total payments from all prosumers exceed the total community cost, which is captured by the energy community manager’s objective function. 

Note that the bilinear term $x_{n,t} \ y_{n,t}$, representing the payment of prosumer $n$ at time $t$, appears in both \eqref{eq:upper_IR} and \eqref{eq:upper_RA}. To avoid nonlinearity, this bilinear term will be replaced with an exact linear reformulation, which is derived using strong duality theory in the lower-level. This reformulation is presented in Section \ref{subsubsec:reformulation} 

\subsubsection{The lower-level problems: Prosumers' power signature} \label{subsubsec:consumer}

Each prosumer's load is the weighted sum of an unknown set of possible signatures. Alternatively, one can say that each prosumer's load is a weighted sum of all signatures in the set $\mathcal{K}$, where many weights will turn out to be zero. In this subsection, we will provide a mathematical model of all these possible signatures, which will each be modeled as one lower-level problem. Each lower-level problem aims to minimize the electricity procurement cost for a specific possible signature $k$ of prosumer $n$. The objective function of each lower-level problem is the product of the price $x_{n,t}$ for prosumer $n$ at time  $t$ and its corresponding power import/export $p_{n,t,k}$ associated with signature $k$. The response of prosumer $n$ for signature $k$ is therefore the optimal solution to the following lower-level optimization problem:
\begin{subequations}
\label{eqs:lower_level_general}
    \begin{align}
    \underset{p_{n,t,k}}{\min} \: & \sum_{t\in\mathcal{T}} x_{n,t} p_{n,t,k}\\
    \mathrm{s.t.} \: & \ \ \{p_{n,t,k}\}_{t \in \mathcal{T}} \in \mathcal{F}_{n,k} \left(\Gamma_{n,k}\right),
    \end{align}
\end{subequations}
where the feasible space $\mathcal{F}_{n,k}\left(\Gamma_{n,k}\right)$ represents a set of constraints defining the feasible power import/export for signature $k$ of prosumer $n$, with $\Gamma_{n,k}$ denoting the input parameters that characterize the technical characteristics and usage preferences of the prosumer's assets. Below, we detail the expression of these constraints and input parameters for five signature types, representing the most common types of flexible loads and DERs. Additionally, we list the dual variables associated with each constraint, which will be used to formulate the Karush-Kuhn-Tucker (KKT) conditions for each lower-level problem. 

\textit{Flexible baseload signature type}: A signature $k$ of this type represents the base consumption of a prosumer $n$, typically a residential user, including lighting, kitchen appliances, and other necessities. These assets are generally inflexible throughout the day, except for a short window of time $\mathcal{T}_{n,k}^{{\rm{flex}}}$ when the prosumer may manually shift this load. The constraints governing this signature model the feasible imported/exported power, while respecting flexibility and operational limits on the base consumption.
The following constraints define the feasible space $\mathcal{F}_{n,k}\left(\Gamma_{n,k}\right)$ for this signature type. 
Constraint \eqref{eq:flex_balance} ensures that the imported power associated with this signature matches the optimal base load profile $l_{n,t,k}$ throughout the day:
\begin{subequations}
\begin{align}
& - p_{n,t,k} + l_{n,t,k}= 0 && \forall t \in \mathcal{T} && :\lambda^{(1)}_{n,t,k}. \label{eq:flex_balance}
\end{align}

Constraint \eqref{eq:nonflex_hours} ensures that the optimal load profile during non-flexible hours remains equal to the unshifted base load prediction $\mathrm{L}_{n,t}$, which is an input parameter:
\begin{align}
& l_{n,t,k} - \mathrm{L}_{n,t} = 0 && \forall t \in \mathcal{T} \setminus \mathcal{T}_{n,k}^{{\rm{flex}}} && \ :\lambda^{(2)}_{n,t,k}. \label{eq:nonflex_hours}
\end{align}

In contrast,  \eqref{eq:flex_hours} ensures that the optimal load consumed during the flexible hours matches the unshifted base load prediction for that time period:

\begin{align}
& \sum_{t \in \mathcal{T}_{n,k}^{{\rm{flex}}}} \left( l_{n,t,k} - \mathrm{L}_{n,t} \right) = 0 && && :\lambda^{(3)}_{n,k}. \label{eq:flex_hours}
\end{align}

Lastly,  \eqref{eq:flex_bounds} ensures that the optimal shifted load at each time period remains within the bounds of the maximum predicted load $\overline{\mathrm{L}}_n$ and the minimum predicted load $\underline{\mathrm{L}}_n$:
\begin{align}
& \underline{\mathrm{L}}_n \leq l_{n,t,k} \leq \overline{\mathrm{L}}_n && \forall t \in \mathcal{T} && :\underline{\mu}^{{\rm{L}}}_{n,t,k}, \overline{\mu}^{{\rm{L}}}_{n,t,k}. \label{eq:flex_bounds}
\end{align}
\end{subequations}

Constraints \eqref{eq:flex_balance}–\eqref{eq:flex_bounds} collectively define the set of feasible power import/export for this signature type as:
\begin{equation}
\label{eqs:flexible hours set}
\mathcal{F}_{n,k}^{\rm{L}} \left(\underline{\mathrm{L}}_n,\overline{\mathrm{L}}_n,\{\mathrm{L}_{n,t}\}_{t \in \mathcal{T}},\mathcal{T}_{n,k}^{{\rm{flex}}}\right) = \left\{ \{p_{n,t,k}\}_{t \in \mathcal{T}} \in \mathbb{R}^{|\mathcal{T}|} \mid \exists \ \{l_{n,t,k}\}_{t \in \mathcal{T}} \in \mathbb{R}^{|\mathcal{T}|} \ {\rm{s.t.}} \ \eqref{eq:flex_balance}-\eqref{eq:flex_bounds} \right\}.
\end{equation}

\textit{PV signature type}: A signature $k$ of this type models the production of a PV system owned by a prosumer $n$. There is no inherent flexibility in this signature, as curtailing PV power production goes against the principles of an energy community and is not generally preferred. Therefore, the only constraint associated with this signature ensures that the hourly PV power production profile $\mathrm{PV}_{n,t}$, which is an input parameter, is equal to the variable power import/export of the prosumer:
%
\begin{align}
& - p_{n,t,k} - \mathrm{PV}_{n,t} = 0 \ &&\forall t \in \mathcal{T} &&:\lambda^{(4)}_{n,t,k}. \label{eq:solar}
\end{align}

This constraint defines the set of feasible power import/export for this signature type as:
\begin{equation}
\label{eqs:pv}
\mathcal{F}_{n,k}^{\mathrm{PV}} \left(\{\mathrm{PV}_{n,t}\}_{t \in \mathcal{T}}\right) = \left\{ \{p_{n,t,k}\}_{t \in \mathcal{T}} \in \mathbb{R}^{|\mathcal{T}|} \mid \eqref{eq:solar} \right\}.
\end{equation}

The reference hourly PV production profile, $\mathrm{PV}_{n,t}$, for each prosumer is defined as an input parameter of this signature type, so the PV capacity does not need to be specified. In the case study, we will further discuss how varying PV capacities are handled across the community.

\textit{Battery signature type}: A signature $k$ of this type models the optimal usage of a home battery system by prosumer $n$. The constraints for this signature define the feasible imported/exported power while respecting the operational limits of the battery.
Constraint \eqref{eq:bat_balance} ensures that the import/export from the grid is equal to the battery charging or discharging power, denoted as $b_{n,t,k}$:
\begin{subequations}
\begin{align}
& - p_{n,t,k} + b_{n,t,k} = 0 \ && \forall t \in \mathcal{T} && :\lambda^{(5)}_{n,t,k}. \label{eq:bat_balance}
\end{align}

Constraints \eqref{eq:bat_soe1}-\eqref{eq:bat_soe3} track the state of energy  $e_{n,t,k}$ of the battery, ensuring that the final state of energy is equal to the initial state, $E^{0}_n$, which is an input parameter. This condition allows for continuity of the battery's state across days:
\begin{align}
& e_{n,t-1,k} - e_{n,t,k} + b_{n,t,k}= 0 && \forall t \in \mathcal{T} \setminus \{1\} &&:\lambda^{(6)}_{n,t,k} \label{eq:bat_soe1}\\
& E^{0}_n - e_{n,t=1,k} + b_{n,t=1,k} = 0 && && :\lambda^{(6)}_{n,t=1,k} \label{eq:bat_soe2}\\
& e_{n,t=24,k} - E^{0}_n = 0 && && :\lambda^{(7)}_{n,k}. \label{eq:bat_soe3}
\end{align}

Finally, \eqref{eq:bat_limits1} and \eqref{eq:bat_limits2} ensure that the battery's discharging and charging limits  ($\underline{\mathrm{B}}_n,\overline{\mathrm{B}}_n$), as well as the maximum and minimum storage capacity limits ($\underline{\mathrm{E}}_n, \overline{\mathrm{E}}_n$), which are input parameters, are not violated:
\begin{align}
    & \underline{\mathrm{B}}_n \leq b_{n,t,k} \leq \overline{\mathrm{B}}_n && \forall t \in \mathcal{T} && :\underline{\mu}^{{\rm{B}}}_{n,t,k}, \overline{\mu}^{{\rm{B}}}_{n,t,k} \label{eq:bat_limits1} \\
    & \underline{\mathrm{E}}_n \leq e_{n,t,k} \leq \overline{\mathrm{E}}_n && \forall t \in \mathcal{T} && :\underline{\mu}^{{\rm{E}}}_{n,t,k}, \overline{\mu}^{{\rm{E}}}_{n,t,k}.\label{eq:bat_limits2}
\end{align}
\end{subequations}
These constraints collectively define the set of feasible power import/export for this signature type as:
\begin{equation}
\label{eqs:battery}
\begin{split}
& \mathcal{F}_{n,k}^{\rm{E}} \left( E^{0}_n, \underline{\mathrm{B}}_n,\overline{\mathrm{B}}_n , \underline{\mathrm{E}}_n, \overline{\mathrm{E}}_n \right) \\
& = \left\{ \{p_{n,t,k}\}_{t \in \mathcal{T}} \in \mathbb{R}^{|\mathcal{T}|} \mid \exists \ \{b_{n,t,k}\}_{t \in \mathcal{T}} \in \mathbb{R}^{|\mathcal{T}|} , \{e_{n,t,k}\}_{t \in \mathcal{T}} \in \mathbb{R}^{|\mathcal{T}|} \ {\rm{ s.t.}} \ \eqref{eq:bat_balance}-\eqref{eq:bat_limits2} \right\}.
\end{split}
\end{equation}

\textit{Heat pump signature type}: A signature $k$ of this type models the flexible power consumption of a heat pump used for heating the residence of a prosumer $n$. The constraints modeling this signature type are inspired by the work of \cite{Hao2015TCL}, which employs a generalized battery model to formulate a linear program for thermostatically controlled loads.
Constraint \eqref{eq:tcl_balance} ensures that the grid power import/export is equal to the power consumed by the heat pump, denoted as $p^{\mathrm{TCL}}_{n,t,k}$:
\begin{subequations}
\begin{align}
& - p_{n,t,k}  + p^{\mathrm{TCL}}_{n,t,k}= 0 \ && \forall t \in \mathcal{T} && :\lambda^{(8)}_{n,t,k}. \label{eq:tcl_balance}
\end{align}

Constraints \eqref{eq:tcl_temp1}-\eqref{eq:tcl_temp3} track the hourly indoor temperatures $\tau_{n,t,k}$ and account for heat loss or gain due to the temperature difference between the indoor and ambient temperatures $\tau^{\mathrm{ex}}_t$. These constraints are based on initial and final indoor temperature parameters $\tau^{0}_n$ and $\tau^{\mathrm{end}}_n$, as well as the heat pump parameters, such as the coefficient of performance $\eta_n$, and the building’s thermal parameters, including thermal resistance $R_n$ and thermal capacity $C_n$:
\begin{align}
    & \tau_{n,t,k} = \tau_{n,t-1,k} - \frac{1}{R_n C_n} \left( \tau_{n,t-1,k} - \tau^{\mathrm{ex}}_t  \right) + \frac{\eta_n}{C_n} p^{\mathrm{TCL}}_{n,t,k} && \forall t \in \mathcal{T} \setminus 1 && :\lambda^{(9)}_{n,t,k} \label{eq:tcl_temp1}\\
    & \tau_{n,t=1,k} = \tau_n^0 - \frac{1}{R_n C_n} \left( \tau^{0}_n - \tau^{\mathrm{ex}}_1  \right) + \frac{\eta_n}{C_n} p^{\mathrm{TCL}}_{n,t=1,k} && && :\lambda^{(9)}_{n,t=1,k} \label{eq:tcl_temp2} \\
    & \tau_{n,t=24,k} = \tau^{0}_n  && && :\lambda^{(10)}_{n,k}. \label{eq:tcl_temp3}
\end{align}

Constraint \eqref{eq:tcl_power_limit} limits the power consumption of the heat pump to its  capacity $\overline{\mathrm{P}}^{\mathrm{TCL}}_{n}$, which is an input parameter:
\begin{align}
    & 0 \leq p^{\mathrm{TCL}}_{n,t,k} \leq \overline{\mathrm{P}}^{\mathrm{TCL}}_{n} && \forall t \in \mathcal{T} && :\underline{\mu}^{\mathrm{TCL}},\overline{\mu}^{\mathrm{TCL}}. \label{eq:tcl_power_limit}
\end{align}

Lastly, \eqref{eq:tcl_temp_limit} enforces the preferred minimum and maximum indoor temperatures $(\underline{\tau}_{n,t,k},\overline{\tau}_{n,t,k})$ for the prosumer. Since the heat pump does not have cooling capacity in this model, the indoor temperature is constrained to be lower than the maximum of the peak external temperature $\overline{\tau}^{\mathrm{ex}}$ and the prosumer's preferred maximum temperature, to account for abnormally hot days when the house cannot stay below the preferred temperature:
\begin{align}
    & \underline{\tau}_{n} \leq \tau_{n,t,k} \leq \max \{ \overline{\tau}_{n}, \overline{\tau}^{\mathrm{ex}}_t\} && \forall t \in \mathcal{T} && : \underline{\mu}^{\tau},\overline{\mu}^{\tau}. \label{eq:tcl_temp_limit}
\end{align}
\end{subequations}

These constraints collectively define the set of feasible power import/export for this signature type as:
\begin{equation}
\label{eqs:tcl}
\begin{split}
& \mathcal{F}_{n,k}^{\mathrm{TCL}} \left(\eta_n,R_n,C_n,\tau^0_n,\tau^{\mathrm{end}}_n,\overline{\mathrm{P}}^{\mathrm{TCL}}_{n},\{\underline{\tau}_{n,t,k},\overline{\tau}_{n,t,k}\}_{t \in \mathcal{T}}\right) \\
& = \left\{ \{p_{n,t,k}\}_{t \in \mathcal{T}} \in \mathbb{R}^{|\mathcal{T}|} \mid \exists \ \{p^{\mathrm{TCL}}_{n,t,k}\}_{t \in \mathcal{T}} \in \mathbb{R}^{|\mathcal{T}|}, \{\tau_{n,t,k}\}_{t \in \mathcal{T}} \in \mathbb{R}^{|\mathcal{T}|} \ {\rm{ s.t.}} \ \eqref{eq:tcl_balance}-\eqref{eq:tcl_temp_limit} \right\}.    
\end{split}
\end{equation}

\textit{EV signature type}: A signature $k$ of this type models the power consumption of an EV owned by prosumer $n$. The constraints governing this signature model the charging and discharging power of the EV's battery, which is assumed to be capable of delivering power to the grid when plugged in. In contrast to the home battery signature type, these constraints must account for the presence of the electric vehicle at the charging station, represented by a binary input parameter $\mathrm{U}^{\mathrm{EV}}_{n,t,k}$, which is equal to 1 if the vehicle is plugged in and 0 otherwise, along with its energy usage during each trip.
Constraint \eqref{eq:ev_balance} ensures that the grid import/export matches the vehicle charging or discharging power $ev_{n,t,k}$:
\begin{subequations}
\begin{align}
& - p_{n,t,k} + ev_{n,t,k} = 0 \ &&\forall t \in \mathcal{T} && :\lambda^{(11)}_{n,t,k}. \label{eq:ev_balance}
\end{align}

Constraints \eqref{eq:ev_soe1}-\eqref{eq:ev_soe3} track the state of energy of the vehicle battery $s_{n,t,k}$, depending on its charging or discharging power $ev_{n,t,k}$ and its energy usage, computed as the average hourly power usage while driving $\mathrm{P}^{\mathrm{DR}}_n$ multiplied by the binary parameter representing the vehicle absence from its charging station $(1-\mathrm{U}^{\mathrm{EV}}_{n,t,k})$. Adding this term ensures that the EV returns with a lower state of charge after each trip. Additionally, the state of energy of the EV at the end of the day is set to its initial value $S_n^0$ to ensure driving needs can be met the following day:
\begin{align}
    & s_{n,t,k} = s_{n,t-1,k} + ev_{n,t,k} - \mathrm{P}_n^{\mathrm{DR}}\big(1-\mathrm{U}^{\mathrm{EV}}_{n,t,k}\big) &&\forall t \in \mathcal{T} \setminus 1 &&:\lambda^{(12)}_{n,t,k} \label{eq:ev_soe1}\\
    & s_{n,t=1,k} = S_n^0 + ev_{n,t=1,k} - \mathrm{P}^{\mathrm{DR}}_n\big(1-\mathrm{U}^{\mathrm{EV}}_{n,t=1,k}\big) && &&:\lambda^{(12)}_{n,1,k} \label{eq:ev_soe2}\\
    & s_{n,t=24,k} - S_n^0  = 0 && &&:\lambda^{(13)}_{n,k}. \label{eq:ev_soe3}
\end{align}

Finally, \eqref{eq:ev_limits1} and \eqref{eq:ev_limits2} ensure that the EV's charging and discharging limits ($\overline{\mathrm{EV}}_n,\underline{\mathrm{EV}}_n$), minimum and maximum battery capacity limits ($\underline{\mathrm{S}}_n, \overline{\mathrm{S}}_n$), are not violated, whether it is connected or not:
\begin{align}
    & \mathrm{U}^{\mathrm{EV}}_{n,t,k} \underline{\mathrm{EV}}_n \leq ev_{n,t,k} \leq \mathrm{U}^{\mathrm{EV}}_{n,t,k} \overline{\mathrm{EV}}_n && \forall t \in \mathcal{T} && :\underline{\mu}^{{\rm{EV}}}_{n,t,k},\overline{\mu}^{{\rm{EV}}}_{n,t,k} \label{eq:ev_limits1} \\
    & \underline{\mathrm{S}}_n \leq s_{n,t,k} \leq \overline{\mathrm{S}}_n && \forall t \in \mathcal{T} && :\underline{\mu}^{{\rm{S}}}_{n,t,k}, \overline{\mu}^{{\rm{S}}}_{n,t,k}. \label{eq:ev_limits2}
\end{align}
\end{subequations}

These constraints collectively define the set of feasible power import/export for this signature type as:
\begin{equation}
\label{eqs:ev}
\begin{split}
& \mathcal{F}_{n,k}^{\mathrm{EV}} \left(\underline{\mathrm{EV}}_n,\overline{\mathrm{EV}}_n,\underline{\mathrm{S}}_n,\overline{\mathrm{S}}_n,S^0_n\{\mathrm{U}^{\mathrm{EV}}_{n,t,k}\}_{t \in \mathcal{T}}\right) \\
& = \left\{ \{p_{n,t,k}\}_{t \in \mathcal{T}} \in \mathbb{R}^{|\mathcal{T}|} \mid \exists \ \{ev_{n,t,k}\}_{t \in \mathcal{T}} \in \mathbb{R}^{|\mathcal{T}|}, \{s_{n,t,k}\}_{t \in \mathcal{T}} \in \mathbb{R}^{|\mathcal{T}|} \ {\rm{ s.t.}} \ \eqref{eq:ev_balance}-\eqref{eq:ev_limits2} \right\}.    
\end{split}
\end{equation}

\subsubsection{Bilevel problem reformulation} \label{subsubsec:reformulation}

Firstly, the two nonlinear constraints \eqref{eq:upper_IR} and \eqref{eq:upper_RA} in the upper-level problem are reformulated using the principle of strong duality theory. The energy cost of prosumers, which is the primal objective function of the lower-level cost-minimization problem, contains bilinear terms, as it multiplies the price $x_{n,t}$ by the prosumer response $y_{n,t}$. However, due to strong duality theory, these bilinear terms can be replaced with the dual objective functions, which are linear in our cases. This equality ensures that the optimal primal and dual objective functions yield the same value. The bilinear primal objective function used in \eqref{eq:upper_IR} and \eqref{eq:upper_RA} and the equivalent combination of the linear dual objective functions is presented in \eqref{eq:strong_duality}: 

\begin{align} \label{eq:strong_duality}
    \nonumber \sum_{t \in \mathcal{T}} &x_{n,t} y_{n,t} \\ 
    \nonumber =& - \sum_{k\in\mathcal{K}^{\mathrm{L}}}\theta_k \left[\sum_{t \in \mathcal{T} \setminus \mathcal{T}_k^{\mathrm{flex}}} \lambda^{(2)}_{n,t,k} \text{L}_{n,t} + \lambda^{(3)}_{n,k} \sum_{t \in T^{{\mathrm{flex}}}_k} \text{L}_{n,t} - \sum_{t \in \mathcal{T}} \left(\underaccent{\bar}{\mu}^{{\mathrm{L}}}_{n,t,k} \underaccent{\bar}{\text{L}}_n + \bar{\mu}^{{\mathrm{L}}}_{n,t,k} \bar{\text{L}}_n \right)\right] \\
    \nonumber &- \sum_{k \in \mathcal{K}^{\mathrm{PV}}} \theta_{k} \sum_{t \in \mathcal{T}}\lambda^{(4)}_{n,t,k} \mathrm{PV}_{n,t} \\
    \nonumber &+\sum_{k \in \mathcal{K}^{\mathrm{B}}} \theta_k \left[ \lambda^{(6)}_{n,t=1,k} \mathrm{E}^0_n - \lambda^{(7)}_{n,k} \mathrm{E}^0_n + \sum_{t \in \mathcal{T}} \left(\underline{\mu}^{{\mathrm{B}}}_{n,t,k} \underline{\mathrm{B}} - \overline{\mu}^{{\mathrm{B}}}_{n,t,k} \overline{\mathrm{B}} + \underline{\mu}^{{\mathrm{E}}}_{n,t,k} \underline{\mathrm{E}} - \overline{\mu}^{{\mathrm{E}}}_{n,t,k} \overline{\mathrm{E}}\right)\right]\\
    \nonumber &+ \sum_{k \in \mathcal{K}^{\mathrm{TCL}}} \theta_k \Bigg[ \lambda^{(9)}_{n,t=1,k} \left( \tau^0_n - \frac{\tau^0_n}{R_n C_n}\right) - \lambda^{(10)}_{n,k} \tau^0_n + \sum_{t \in \mathcal{T}} \Big( \lambda^{(9)}_{n,t,k} \frac{\tau^{\mathrm{ex}}_t}{R_n C_n} + \underline{\mu}^{\tau} \underline{\tau}_{n,k} \\
    \nonumber & \hspace{240pt} - \overline{\mu}^{\tau} \overline{\tau}_{n,k} - \overline{\mu}^{\mathrm{TCL}} \overline{\mathrm{P}}^{\mathrm{TCL}}_n \Big) \Bigg]\\
    \nonumber &- \sum_{k \in \mathcal{K}^{\mathrm{EV}}} \theta_k \Bigg[\ \mathrm{S}^0_n \left( \lambda^{(13)}_{n,k} + \lambda^{(12)}_{n,t=1,k}\right) + \sum_{t \in \mathcal{T}} \Big(-\lambda^{(12)}_{n,t,k} \mathrm{P}^{\mathrm{DR}}_n (1 - \mathrm{U}_{n,t,k}^{\mathrm{EV}}) + \underline{\mu}^{{\mathrm{EV}}}_{n,t,k} \mathrm{U}_{n,t,k}^{\mathrm{EV}} \underline{\mathrm{EV}}_n \\
    & \hspace{215pt}- \overline{\mu}^{{\rm{EV}}}_{n,t,k} \mathrm{U}^{\mathrm{EV}}_{n,t,k} \overline{\mathrm{EV}}_n + \underline{\mu}^{{\rm{S}}}_{n,t,k} \underline{\mathrm{S}}_n - \overline{\mu}^{{\rm{S}}}_{n,t,k} \overline{\mathrm{S}}_n \Big) \Bigg].
\end{align}

To reformulate the bilevel problem as a single-level optimization, we apply the KKT conditions to all lower-level optimization problems and incorporate them into the upper-level problem. Since the lower-level problems are linear programs, the KKT conditions serve as both necessary and sufficient optimality conditions. Incorporating these conditions as constraints in the upper-level problem transforms it into a mathematical program with equilibrium constraints (MPEC). Using the well-known Fortuny-Amat approach, this MPEC is reformulated into a mixed-integer linear program (MILP), which can be solved in a single step and seamlessly integrated into a learning framework. We refer to the reformulated MILP as \textbf{BiPS}, short for Bilevel Price Setting.

\subsection{Learning algorithm}

The general framework for Thompson sampling was introduced in Section \ref{sec:framework}. This section aims to provide detailed mathematical formulations of the different steps of the implemented algorithm, as outlined in Algorithm \ref{alg:TS}. Throughout the algorithm definition, we refer back to Figure \ref{fig:ts-framework} to establish the link between the proposed algorithm and the generalized Thompson sampling framework.

The initial step (\textbf{Step 1}) of the algorithm requires collecting all available information about known parameters in the community, such as DER capacities. The algorithm does not rely on this information to perform well. Still, leveraging any prior information could improve the speed of learning and minimize regret in early iterations. In this initial step, the prior distributions for the unknown parameters $f_d(\boldsymbol{\theta}_n)$ are initialized. For easy updates, a Gaussian normal distribution is initialized for all unknown weights, i.e., one for each signature and prosumer. As the true values for the weights are generally assumed to lie between 0 and 1, the distributions are initialized such that most generated samples will lie within this range.

\begin{algorithm}[t]
    \caption{Thompson sampling algorithm for learning prosumer behavior}
    \begin{algorithmic}
        \Input Prosumer parameter set $\boldsymbol{\Omega}$
        \Initialize \For{$n = 1,...,N$} 
            \State $f_{d=1}(\boldsymbol{\theta}_n) \sim \mathcal{N}(\boldsymbol{\overline{\theta}}_{n},\boldsymbol{\Sigma}_n)$ (\textbf{Step 1} in Figure \ref{fig:ts-framework}) 
        \EndFor
        \For{$d = 1, ..., \mathcal{D}$}
        \State $\mathbf{for} \ n = 1,...,N$ \\
            \hspace{30pt}1. Sample daily weights $\tilde{\boldsymbol{\theta}}_{n,d}$ from the prior distributions $f_{d}(\boldsymbol{\theta}_n)$ (\textbf{Step 2} in Figure \ref{fig:ts-framework})
            \State \textbf{end for}
            
            \State 2. Observe daily exogenous parameter set $\boldsymbol{\Phi}_d$
            
            \State 3. Determine prices according to sampled values $\tilde{\boldsymbol{\Theta}}_d$ (\textbf{Step 3} in Figure \ref{fig:ts-framework}):
            \begin{align*}
                \tilde{\boldsymbol{X}}_d = \underset{\textbf{X}}{\mathrm{argmin}} \{ \textbf{BiPS} \left(\mathbf{\Omega},\mathbf{\Phi}_d,\boldsymbol{\theta} \right) \mid \boldsymbol{\theta} = \tilde{\boldsymbol{\Theta}}_d \}
            \end{align*}
            
            \State $\mathbf{for} \ n = 1, ..., N$ \\
            \hspace{30pt}4. Observe response of prosumers to prices resulting from sample $\tilde{\boldsymbol{X}}_d$ (\textbf{Step 4} in Figure \ref{fig:ts-framework}): 
            \begin{equation*}
                \tilde{\boldsymbol{y}}_{n,d} = \tilde{\mathbf{P}}_n \boldsymbol{\theta}_n^* + \boldsymbol{\epsilon} \quad \mathrm{where} \quad \boldsymbol{\epsilon} \sim \mathcal{N}(0,\sigma)
            \end{equation*} \\
            \State \textbf{end for}

            \State $\mathbf{for} \ n = 1, ..., N$ \\
            \hspace{30pt}5. Update the mean and covariance of the posterior distribution $f_{d+1}(\boldsymbol{\theta}_n)$(\textbf{Step 5} in Figure \ref{fig:ts-framework})
            \begin{equation} \label{eq:mean_update}
            \boldsymbol{\overline{\theta}}_{n,d+1} = \boldsymbol{\overline{\theta}}_{n,d} + \boldsymbol{\Sigma}_{n,d} \, \tilde{\mathbf{P}}_n^\intercal (\tilde{\mathbf{P}}_n \boldsymbol{\Sigma}_{n,d} \, \tilde{\mathbf{P}}_n^\intercal + \textbf{R}_n)^{-1} (\tilde{\boldsymbol{y}}_{n,d} - \tilde{\mathbf{P}}_n \boldsymbol{\overline{\theta}}_{n,d})
            \end{equation}
            \begin{equation} \label{eq:covar_update}
            \boldsymbol{\Sigma}_{n,d+1} = \boldsymbol{\Sigma}_{n,d} - \boldsymbol{\Sigma}_{n,d} \, \tilde{\mathbf{P}}_n^\intercal (\tilde{\mathbf{P}}_n \boldsymbol{\Sigma}_{n,d}\, \tilde{\mathbf{P}}_n^\intercal + \textbf{R}_n)^{-1} \tilde{\mathbf{P}}_n \boldsymbol{\Sigma}_{n,d}
            \end{equation}
            \State \textbf{end for}
        \EndFor \ (\textbf{Step 6} in Figure \ref{fig:ts-framework})
        
    \end{algorithmic}
    \label{alg:TS}
\end{algorithm}

Once the prior distribution is initialized, the iterative steps (\textbf{Step 2 - Step 6}) of the Thompson sampling algorithm are initiated. In \textbf{Step 2}, a daily sample of all unknown parameters is taken from their prior distribution. Additionally, all daily exogenous parameters, such as PV power production, outdoor temperature, and day-ahead prices, are observed. In \textbf{Step 3}, the  community manager can now set the optimal dynamic prices while anticipating the responses of the prosumers by solving the bilevel problem \textbf{BiPS} with these sampled and exogenous parameters as inputs. These prices $\tilde{\textbf{x}}_d$ are then shared with the prosumers. Following this, in \textbf{Step 4}, the true response from all prosumers is observed through their individual smart meters. This response is the weighted sum of the power signatures resulting from the different lower-level optimization problems with some additional Gaussian white noise $\boldsymbol{\epsilon}$ as explained in equation \eqref{eq:response_model}. It is important to note that the weight used in the observation calculation is not the same as the samples used in the optimization, but rather the true underlying weights that are unknown to the community manager. Once the observations for the daily prices have been gathered for all hours of the day, the posterior distributions of the unknown parameters are updated in \textbf{Step 5}, using Bayesian linear regression. The mathematical formulation of the mean and covariance updates is provided in equations \eqref{eq:mean_update} and \eqref{eq:covar_update}. Once the update has taken place, the algorithm progresses to the next iteration, in \textbf{Step 6}, where the updated posterior distribution becomes the prior and the iterative process starts again with taking a sample. Over iterations, the posterior distributions should converge to the true value of the unknown parameters.

\subsection{Regret metric}

To accurately assess the performance of the above learning algorithm over iterations, we track the so-called regret. Regret is the cumulative cost (or reward) difference between the implemented learning algorithm and an optimal algorithm with perfect knowledge of the true weights of all signatures for all prosumers. This optimal algorithm is able to make the best possible choices at all times. Regret is not a metric that is tracked by the learner that implements an algorithm, as the cost or reward of the optimal algorithm is not accessible to the learner. While some regret metrics do consider noise, the regret metric in this work does not, and the expected responses are used to calculate the regret.

The regret $\mathcal{R}_{d}$ in iteration $d$ considers the difference in community cost (the upper-level objective function) under Thompson sampling versus under perfect information, as shown in \eqref{eq:regret}. For brevity, we refer to the total community cost as $\mathrm{cost}^{\rm{C}}({\cdot})$ instead of writing out the upper-level objective function \eqref{eq:upper_obj}:

\begin{equation} \label{eq:regret}
    \mathcal{R}_{d} = \mathrm{cost}^{\rm{C}} \left(\tilde{\boldsymbol{x}}_{nd} \mid \boldsymbol{\theta}^* \right) - \mathrm{cost}^{\rm{C}} \left( \boldsymbol{x}^*_{nd} \mid \boldsymbol{\theta}^* \right),
\end{equation}
where the first term is the value of the objective function when the prices are optimized with the sampled signature weights, while the prosumers respond with their true signatures $\boldsymbol{\theta}^*$. The second term represents the cost if the true signature weights were known beforehand. As the estimate of the unknown parameters improves, these two terms will converge, meaning that a near-zero regret is obtained in later iterations.

This paper does not discuss or analyze the theoretical bounds of regret under Thompson sampling. Instead, the interested reader is referred to \cite{gopalan} and \cite{Tucker2020ConTS}, which build on each other to provide a thorough discussion of regret bounds in similar learning situations.

\section{Numerical case study} \label{sec:case}

\subsection{Case study setup}

We demonstrate the performance of the proposed algorithm in an energy community with 25 prosumers, each represented by the same set of 10 possible signatures. The signatures are listed below:
\begin{enumerate}[label = ${\arabic*}$:]
    \item Flexible baseload between 6:00 to 10:00 (morning),
    \item Flexible baseload between 10:00 and 17:00 (daytime),
    \item Flexible baseload between 17:00 and 22:00 (evening),
    \item PV power production,
    \item Home battery system,
    \item Heat pump with small acceptable temperature range (19-21$^\circ$C),
    \item Heat pump with large acceptable temperature range (16-24$^\circ$C),
    \item Electric vehicle that is unavailable from 8:00 to 19:00,
    \item Electric vehicle that is unavailable from 6:00 to 15:00 and from 19:00 to 22:00,
    \item Electric vehicle that is unavailable from 7:00 to 10:00 and from 16:00 to 20:00.
\end{enumerate}

The first three signatures are of the flexible baseload signature type. These signatures are differentiated based on the range of flexible hours $\mathcal{T}^{\rm{flex}}_{nk}$ in each signature. For each prosumer $n$, the sum of the true weights over the set of baseload signatures is equal to 1. This ensures that the flexible baseload of each prosumer can be fully represented as a mixture of these signatures, even if they are partially flexible in different hours. In our case study, three different signatures are included as described above.
Signature 4 represents a PV signature type and is characterized by a reference production profile for a 1-kW PV system. The contribution of this signature to each prosumer's power import/export is determined by the true underlying weight, which, in this case, acts as a scaling factor. This scaling factor varies between 0 and 3, reflecting the installed capacity of each prosumer's PV system, accounting for variations in system sizes ranging from 0 kW to 3 kW.
Signature 5 represents a residential battery signature. This asset is unaffected by prosumer behavior, so this case study only includes one standardized battery signature. The true weights associated with this signature can take a value of 0 or 1, representing whether a given prosumer owns a battery or not.
The heat pump and EV signature types are differentiated based on the usage preferences of the prosumers. The heat pump signatures 6 and 7 are differentiated based on two incompatible preferred temperature ranges, represented by the input parameters $\underline{\tau}_{ntk}$ and $\overline{\tau}_{ntk}$. Similarly, the EV signatures 8-10 are differentiated based on three incompatible driving patterns, represented by the input parameters $\mathrm{U}^{\mathrm{EV}}_{ntk}$. Therefore, the true weights associated with these signatures can take a value of 0 or 1, representing the ownership and usage preferences of these assets for each prosumer. If a given prosumer does not own a heat pump (EV) the true weights associated with all the heat pump (EV) signatures will be equal to 0. On the contrary, if a prosumer owns a heat pump (EV) one of the true weights associated with the heat pump (EV) signatures will be equal to 1, representing the prosumer's preferred temperature range (driving pattern). The proposed algorithm therefore aims to learn simultaneously whether a given prosumer owns an asset, and, if so, what their usage preferences are.

\begin{figure}[t]
    \centering
    \input{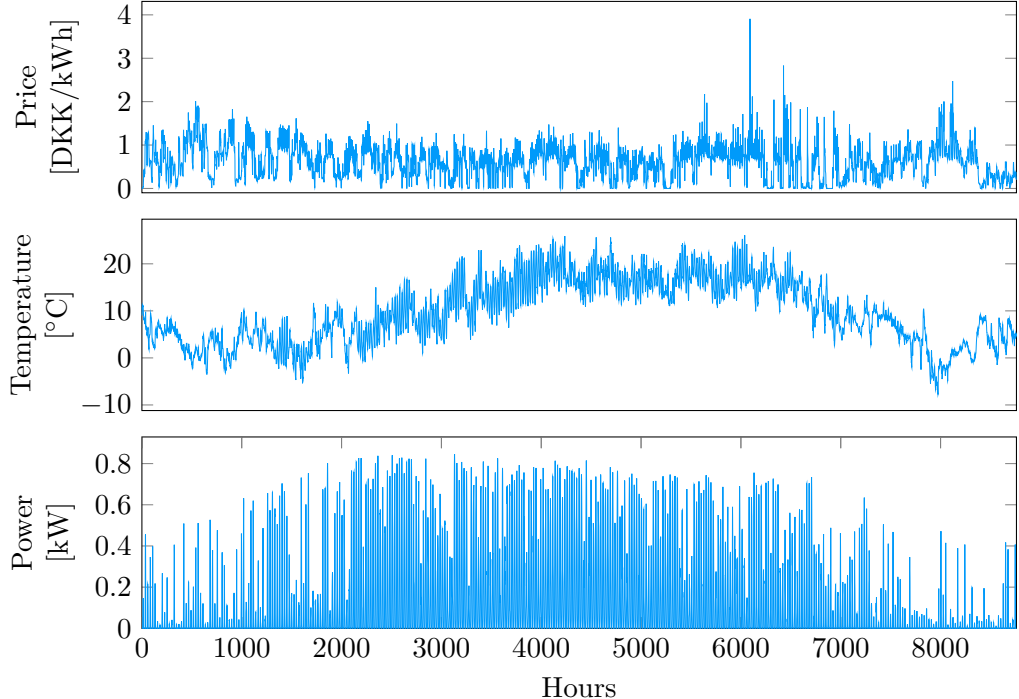}
    \caption{Variations of the exogenous parameter sets $\mathbf{\Phi}_d$ throughout the days of the numerical study.}
    \label{fig:input}
\end{figure}

Furthermore, historical data simulates the exogenous parameters for every day. This means that the algorithm must learn in a context that changes daily. A full year of data is used for day-ahead prices, outdoor temperature, and PV power production. The day-ahead prices have been taken from \cite{Energinet2024Elspot} and filtered to include only the DK2 bidding zone in 2023. Therefore, the prices are in Danish crown per kilowatt-hour [DKK/kWh]. Note that the day-ahead prices have been adjusted not to include negative prices. Weather data is taken from \cite{DMI2024temp} and filtered to be for a location within the same bidding zone and at the same time as the day-ahead prices. The flexible baseload of all prosumers is constant across all days and generated using \cite{McKenna2020}. Lastly, yearly PV power production data is generated using the methodologies described in \cite{Pfenninger2016}. The variation of the exogenous input parameters is illustrated in Figure \ref{fig:input}.

This numerical analysis is run 20 times to ensure the algorithm consistently learns across varying sample sets and mitigates the impact of noise and outliers. The results presented in this section represent the average over these 20 runs. For all runs, the prior distributions of all signature weights are initialized as independent Gaussian distributions $\mathcal{N}(0.5,0.15)$. This ensures that the majority of the prior density function lies within 0 and 1, covering most possible values of the signature weights, while samples outside this range remain rare. The code for this numerical study is provided in an online repository \citep{Crowley_Online_Appendix_to_2025}.

\subsection{Regret}

\begin{figure}[t]
    \centering
    \input{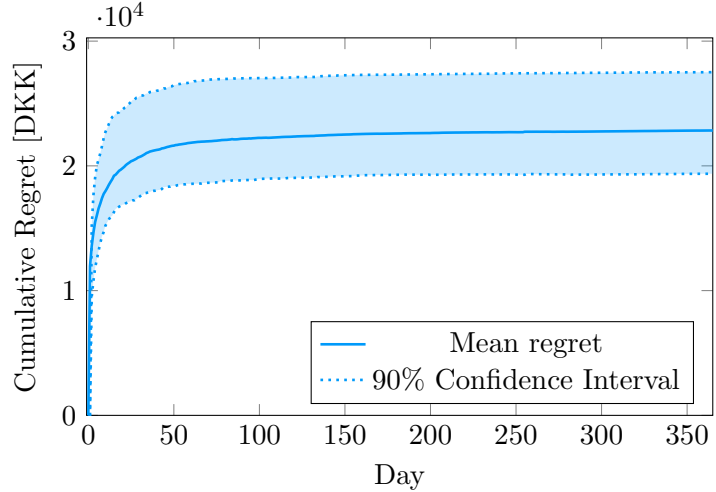}
    \caption{Mean cumulative regret and $90\%$ confidence interval over 365 days}
    \label{fig:regret}
\end{figure}

We first investigate the evolution of the cumulative regret as defined in Section \ref{sec:math}. As illustrated in Figure \ref{fig:regret}, the algorithm consistently achieves near-zero regret after 100 days. Notably, by day 50, the mean cumulative regret curve has already significantly plateaued, with the majority of the learning taking place within the first 25 days. Additionally, Figure \ref{fig:regret} shows the 90\% confidence interval of the cumulative regret across all runs, which exhibits the same behavior beyond 100 days. In fact, the largest difference in regret across any of the runs occurs within the first 25 days. This highlights the impact of the high variance of the prior distribution in the earlier days, which can lead to larger deviations in regret across different runs. As the distributions converge over time, the regret stabilizes across all runs.

\subsection{Convergence of distributions}

\begin{figure}[t]
    \centering
    \input{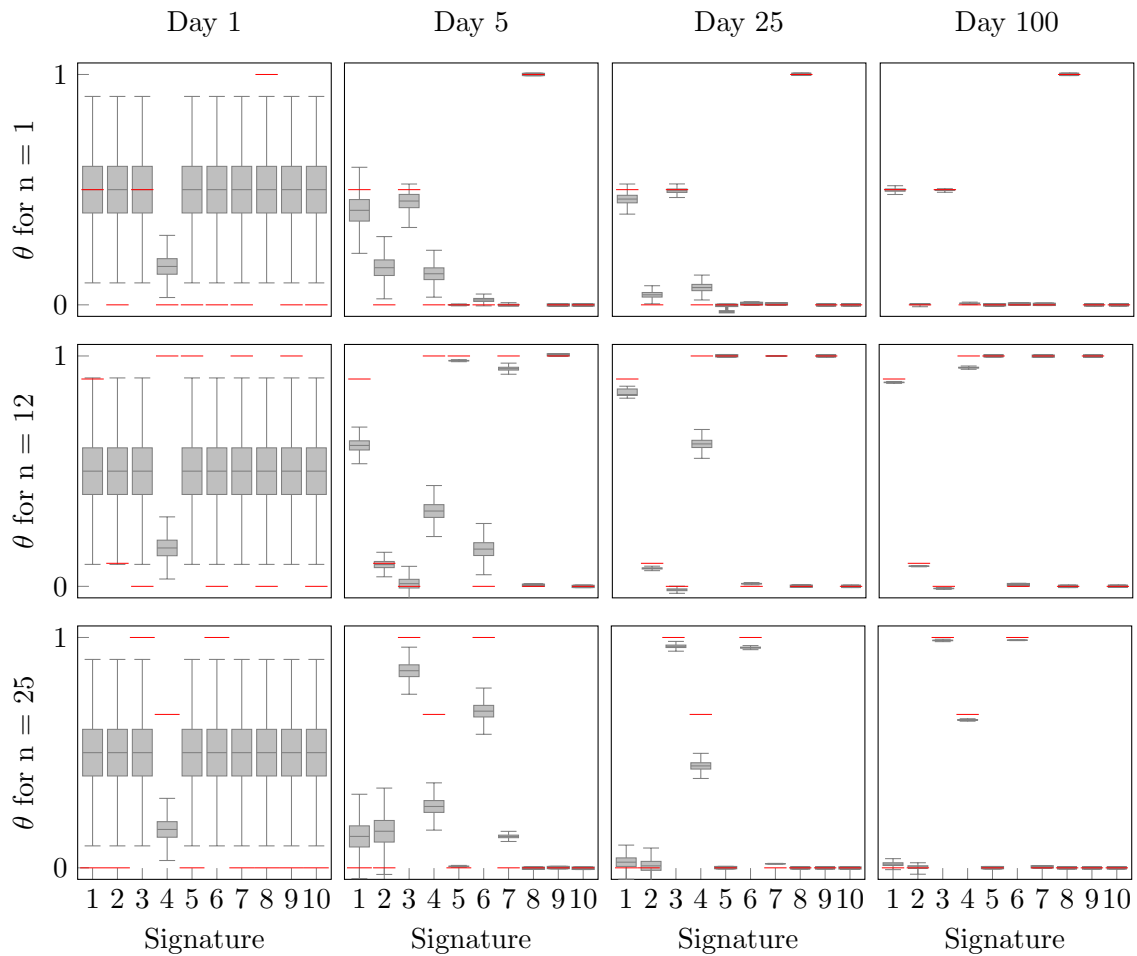}
    \caption{Convergence of the posterior distributions toward the true weights for three representative prosumers over 100 days. Each row of four subplots represents a prosumer. Each column of three subplot represents a specific day. Each box plot approximates the prior/posterior distribution of the weights $\theta$ for ten signatures on a given day and for the given prosumer, while the red line represents the true signature weights $\theta^*$. Note that the weights associated with signature four have been normalized to be between 0 and 1.}
    \label{fig:boxplots}
\end{figure}

The next area investigated is the evolution of the posterior distributions to see if these converge toward the true weights. The convergence of all distributions is illustrated for three representative prosumers across the first 100 days in Figure \ref{fig:boxplots}. One can observe that, on day 1, the prior distributions are the same across all signatures, with the exception of signature 4. This is due to the normalization of the prior distribution of the weights associated with signature 4 to account for the fact that these true weights can be larger than 1. Moving on to day 5, the posterior distributions begin to converge towards the true weights. Specifically, the distributions of the EV and battery signatures have already correctly converged to the true weights, whereas the posterior distributions of all other signatures still have significant variance. By day 25, the heat pump signature weights have also been correctly learned. Finally, by day 100, the posterior distributions of the flexible baseload and PV signature weights have also correctly converged to the true weights. It is worth noting, however, that the posterior distributions for the PV signatures have converged toward a slightly offset mean value for prosumers with a non-zero true weight. We believe this is due to the subpar coverage of the prior, as the algorithm has more difficulty exploring certain regions that are insufficiently covered by the prior distributions. This does, however, show that the algorithm is capable of dealing with a poorly specified prior as long as the prior distribution supports the correct parameter value.

It is worth noting the relationship between the peak power of each signature and the convergence rate of their respective posterior distributions. The charging and discharging powers of the EV and battery signatures have the two highest power contributions, while their distributions also converge the fastest. This trend extends to the heat pump, flexible baseload, and PV signatures, with their peak power and convergence rates following a decreasing order. Learning the true weights associated with the highest power assets faster is beneficial for the community. Indeed, intuitively, by learning the response associated with the higher power assets first, i.e., the EVs and residential batteries, the potential for large violations of the capacity limitation decreases quickly, hereby mitigating the penalties from the DSO, which is the largest potential cost for the energy community. This is further analyzed in the following section.

\subsection{Capacity limit violation}

\begin{figure}[t]
    \centering
    \input{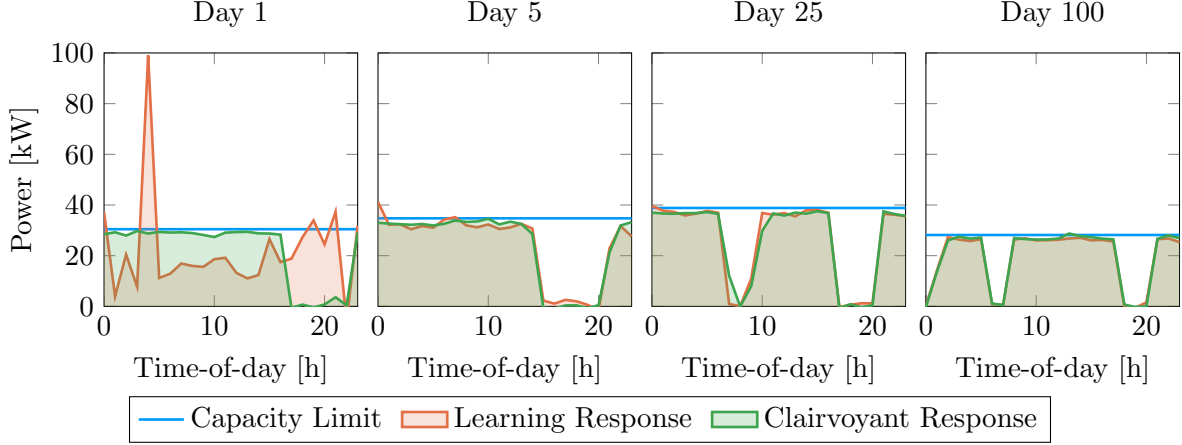}
    \ref{legendlocation}
    \caption{Learned aggregated community response and capacity limitation violation compared to optimal clairvoyant response over the first 100 days.}
    \label{fig:violation}
\end{figure}

We now investigate the grid import/export of the energy community and potential capacity limitation violation, comparing the community's response while learning to the optimal response. As depicted in Figure \ref{fig:violation}, on day 1, the community manager is unable to control the prosumers' response, which greatly differs from the optimal response. This leads to massive peak loads that significantly exceed the capacity limitation in several hours of the day. By day 5, the community's learned response is significantly improved and very similar to the optimal response throughout the day. This results in only a small capacity limitation violation in the first hour of the day. By day 25, the community manager has successfully learned the prosumers' response, optimally controlling it to follow closely the optimal response, with the exception of small deviations in hours 8 and 10. Despite these minor deviations, the capacity limitation is no longer violated throughout the day. Finally, by day 100, even these small differences between the learned and optimal responses have been eliminated, meaning that the learning has succeeded. Figure \ref{fig:violation} reinforces the intuition that the most influential assets with respect to satisfying the capacity limitation are EVs and residential batteries. The community is able to mostly respect the capacity limitation once the true weights of these signatures are correctly learned by day 5. This means that the community manager does not necessaneed to have perfectly learned the true weights of all the flexible assets to mitigate the potentially large penalties resulting from capacity limitation violations and achieve most of the financial benefits for the community. This is a positive sign for practical implementation and widespread acceptance.

\subsection{Non-stationarity}

\begin{figure}[t]
    \centering
    \input{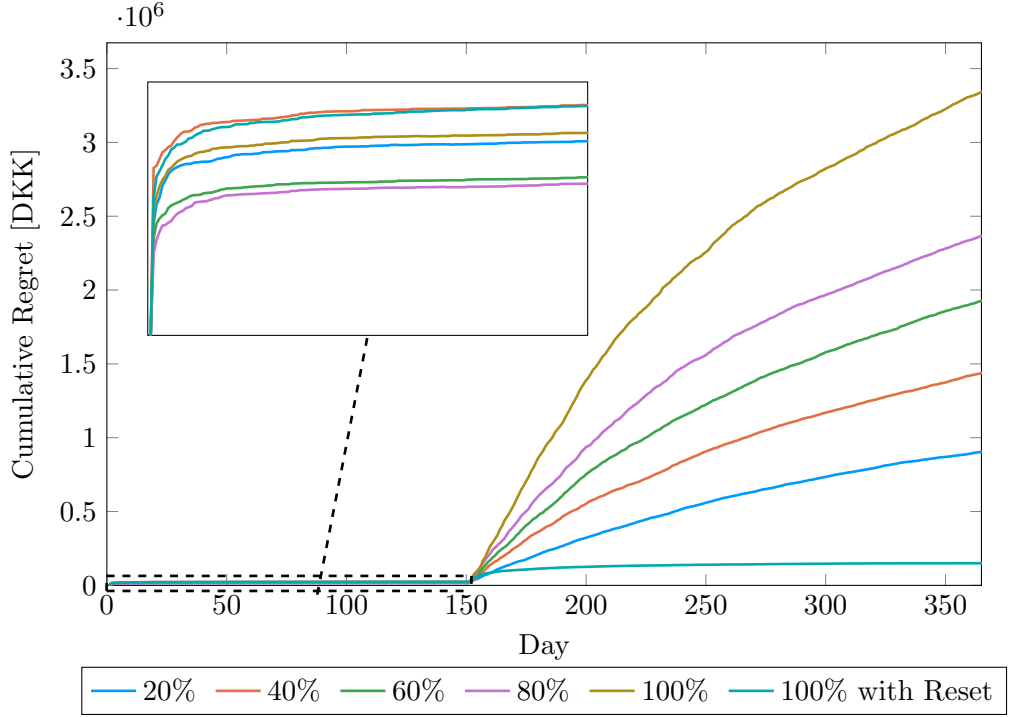}
    \caption{Cumulative regret curve before and after changes in the true weights associated with the heat pump and the EV signatures on day 150.}
    \label{fig:nonstationarity}
\end{figure}

Lastly, we investigate the performance of the proposed algorithm in the case of non-stationarity, i.e., when the true weights associated with the signatures of some of the prosumers change during (or after) the learning phase and the learned weights must be adjusted. To investigate this, five different scenarios are tested in which a varying percentage of the prosumers change their true weights associated with both the heat pump and the EV signatures on day 150, i.e., after the learning has successfully converged and regret has plateaued. As illustrated in Figure \ref{fig:nonstationarity}, in all cases, the regret increases substantially after the non-stationary event occurs. The increase in regret is proportional to the number of prosumers that change their true weight values. Furthermore, the observed regret increase is substantially larger than the regrets achieved with agnostic priors. This is due to the fact that the posterior distributions on day 150 have converged and are generating samples very close to the previous true weights, which are now incorrect. This lack of exploration significantly decreases the learning rate and increases the cumulative regret compared to sampling from a broader prior distribution, which would be more likely to generate samples closer to the new true weights. As a result, we observe that the algorithm cannot attain zero regret again within 215 days after a change in the signature weights occurs. To resolve this, we propose to reset the posterior distributions to the agnostic prior if a significant change in signature weights is detected. This detection could be achieved by looking at the predicted demand profile of a prosumer and comparing it to the observed demand profile. If the distance between these two profiles has increased compared to previous days and is larger than a given noise tolerance, the prior distribution should be reset. Intuitively, the aim of this reset is to reproduce the learning performance of the algorithm that is observed in the first 150 days, as shown in the zoomed window inside Figure \ref{fig:nonstationarity}. Looking at the regret curves in Figure \ref{fig:nonstationarity}, it is clear that resetting the priors is successful and achieves zero regret by approximately day 250, i.e., 100 days after the weights are changed. This corresponds to the same learning time frame observed in Figure \ref{fig:regret}.

\section{Conclusion} \label{sec:conclusion}


This paper introduced a novel application of an online decision-making framework, merging bilevel programming and online learning, specifically Thompson sampling. This novel framework is applied in the context of congestion management in distribution grids, and is used to learn individual prosumers' asset ownership and usage preferences, thereby facilitating the control of the community's grid imports under a price-based demand response scheme. To the authors' knowledge, this is the first framework that integrates bilevel programming into an online learning framework, in addition to learning demand response on an individual level rather than the aggregate level.

The learning performance of this framework is demonstrated in a numerical case study. We show that the implemented algorithm converges to near-zero regret over time. Furthermore, the community is able to learn the true weights associated with high power assets, such as electric vehicles and residential batteries, within 5 days. Doing so enables the community to deliver capacity limitation services to the DSO reliably after less than a week of observations, even before the weights associated with lower power assets have been fully learned. In addition, we observe that after 100 days, even the low power asset weights have been learned correctly. Lastly, the algorithm's learning performance is tested when exposed to a change in the true underlying weights. This simulation showed that a reset of the prior distribution can help achieve near-zero regret again after the true weights have changed. Without this reset, the algorithm's exploration is restricted, impeding its ability to learn effectively.


This work opens several future interesting research directions. Firstly, this numerical analysis is limited to independent normal prior distributions due to their simple updates in linear Bayesian regression. It is of interest to investigate how defining other prior distributions, such as binomial distributions for assets that have weights of either 0 or 1, like heat pumps and electric vehicles, could improve the learning performance of the algorithm in early iterations. Indeed, generating more realistic and feasible samples of the weights using more informed priors in the Thompson sampling algorithm may help reduce regret in early iterations. Secondly, future work should investigate the impact of specifying incorrect sets of signatures in the bilevel model, which do not accurately represent the true assets and usage preferences of the prosumers. In such a case, is the proposed algorithm able to learn the individual responses of the prosumers? Is it capable of mapping these individual responses onto the specified signatures, by finding a different linear combination of weights that closely resembles the unspecified signature? How does this impact the regret and capacity limitation violation of the community? Finally, by identifying the assets of each prosumer and their usage preferences, this work may raise privacy and security concerns. For instance, occupancy information, learned from heat pump or electric vehicle usage, may be sensitive. Therefore, it is interesting to investigate ways to protect the prosumers' sensitive data from such learning algorithms, ensuring that an unwanted third party would be unable to abuse such learned information.

\bibliographystyle{elsarticle-harv} 
\bibliography{bibliography}

\appendix
\section{Nomenclature}
\vspace{-0.75cm}
\nomenclature[A,01]{$\mathcal{D}$}{Set of days (iterations)}
\nomenclature[A,02]{$N$}{Set of  prosumers}
\nomenclature[A,03]{$\mathcal{T}$}{Set of hours}
\nomenclature[A,04]{$\mathcal{T}^{\rm{flex}}_k$}{Subset of flexible hours for signature $k$}
\nomenclature[A,05]{$\mathcal{K}$}{Set of prosumer signatures}
\nomenclature[A,06]{$\mathcal{K}^{(.)}$}{Subset of a specific signature type. The type is specified in the superscript.}
\nomenclature[A,07]{$\mathcal{F}_{nk}^{(.)}$}{Feasible space of  prosumer $n$ for signature $k$}
\nomenclature[B,01]{$\theta_{n,k}$}{Signature weight of prosumer $n$ for signature $k$ [-]}
\nomenclature[B,02]{$\tilde\theta_{n,k}$}{Sample of the signature weight of prosumer $n$ for signature $k$ [-]}
\nomenclature[B,05]{$f_d(\boldsymbol{\theta}_{n})$}{Distribution of signature weights for all signatures of prosumer $n$ on day $d$ [-]}
\nomenclature[B,03]{$\tilde{\boldsymbol{\theta}}_{n}$}{Vector of samples from the prior distribution for prosumer $n$ for all signatures [-]}
\nomenclature[B,04]{$\tilde{\boldsymbol{\Theta}}_{d}$}{Matrix of samples from the prior distribution for all prosumers for all signatures on day $d$ [-]}
\nomenclature[B,06]{$\overline{\boldsymbol{\theta}}_{n,d}$}{Mean vector of the prior distribution for all signatures of prosumer $n$ on day $d$ [-]}
\nomenclature[B,07]{$\boldsymbol{\Sigma}_{n,d}$}{Covariance matrix of the distribution for prosumer $n$ on day $d$ [-]}
\nomenclature[B,08]{$\boldsymbol{\theta}^*_n$}{Vector of true $\theta$ values for all signatures of prosumer $n$ [-]}
\nomenclature[B,09]{$\rho^{\mathrm{spot}}_t$}{Spot price in hour $t$}
\nomenclature[B,10]{$\gamma^{\mathrm{im}}_t$}{Distribution grid import tariff in hour $t$ [DKK/kWh]}
\nomenclature[B,11]{$\gamma^{\mathrm{ex}}_t$}{Distribution grid export tariff in hour $t$ [DKK/kWh]}
\nomenclature[B,12]{$\alpha^{\mathrm{DSO}}_t$}{Penalty for exceeding capacity limitation in hour $t$ [DKK/kW]}
\nomenclature[B,13]{$\overline{\mathrm{P}}^{\mathrm{DSO}}_t$}{Power capacity limit set by the DSO in hour $t$ [kW]}
\nomenclature[B,14]{$\mathrm{cost}^{\mathrm{ex}}_n$}{Energy procurement cost of  prosumer $n$ when it is not part of the community [DKK]}
\nomenclature[B,15]{$\mathrm{L}_{n,t}$}{Base load schedule for prosumer $n$ in hour $t$ [kWh]}
\nomenclature[B,16]{$\mathrm{PV}_{n,t}$}{PV power production for prosumer $n$ in hour $t$ [kWh]}
\nomenclature[B,17]{$\overline{\mathrm{E}}_n$}{Maximum battery state of energy for prosumer $n$ [kWh]}
\nomenclature[B,18]{$\underline{\mathrm{E}}_n$}{Minimum battery state of energy for prosumer $n$ [kWh]}
\nomenclature[B,19]{$\mathrm{E}^0_n$}{Initial battery state of energy for prosumer $n$ [kWh]}
\nomenclature[B,20]{$\overline{\mathrm{B}}_n$}{Maximum battery charging power for prosumer $n$ [kW]}
\nomenclature[B,21]{$\overline{\mathrm{B}}_n$}{Maximum battery discharging power for prosumer $n$ [kW]}
\nomenclature[B,22]{$\tau^{\mathrm{ex}}_t$}{External temperature in hour $t$ [$^{\circ}$C]}
\nomenclature[B,23]{$\underline{\tau}_n$}{Minimum comfortable temperature for prosumer $n$ [$^{\circ}$C]}
\nomenclature[B,24]{$\underline{\tau}_n$}{Maximum comfortable temperature for prosumer $n$ [$^{\circ}$C]}
\nomenclature[B,25]{$\tau^0_n$}{Initial indoor temperature for prosumer $n$ [$^{\circ}$C]}
\nomenclature[B,26]{$\overline{\mathrm{P}}^{\tau}_n$}{Maximum power of the heat pump of prosumer $n$ [kW]}
\nomenclature[B,27]{$\eta^{\tau}_n$}{Heat pump coefficient of performance of prosumer $n$ [-]}
\nomenclature[B,28]{$C_n$}{Thermal capacitance of the residence of prosumer $n$ [kWh/$^{\circ}$C]}
\nomenclature[B,29]{$R_n$}{Thermal resistance of the residence of prosumer $n$ [$^{\circ}$C/kW]}
\nomenclature[B,30]{$\mathrm{P}^{\rm{DR}}_n$}{Average energy usage the EV during an hour of absence for prosumer $n$ [kWh]}
\nomenclature[B,31]{$\overline{\mathrm{S}}_n$}{EV battery energy storage capacity for prosumer $n$[kWh]}
\nomenclature[B,32]{$\underline{\mathrm{S}}_n$}{EV battery minimum energy storage limit for prosumer $n$ [kWh]}
\nomenclature[B,33]{$\mathrm{S}^0_n$}{Initial state of energy in EV battery for prosumer $n$[kWh]}
\nomenclature[B,34]{$\overline{\mathrm{EV}}_n$}{Maximum EV charging power for prosumer $n$ [kW]}
\nomenclature[B,35]{$\underline{\mathrm{EV}}_n$}{Maximum EV discharging power for prosumer $n$ [kW]}
\nomenclature[B,36]{$\mathrm{U}^{\mathrm{EV}}_{n,t,k}$}{Binary parameter for EV connection status for prosumer $n$ in hour $t$ for signature $k$: 1 = connected, 0 = disconnected [-]}
\nomenclature[B,37]{$\boldsymbol{\Omega}$}{Set of Parameters of any prior information for prosumers}
\nomenclature[B,38]{$\boldsymbol{\Phi}_d$}{Exogenous parameters relevant for optimization on day $d$}
\nomenclature[B,39]{$\mathbf{R}_n$}{Covariance matrix of response noise for all signatures of prosumer $n$ [-]}
\nomenclature[C,01]{$x_{n,t}$}{Prices set by the community manager for prosumer $n$ in hour $t$ [DKK/kWh]}
\nomenclature[C,02]{$y_{n,t}$}{Observed response of prosumer $n$ in hour $t$ [kWh]}
\nomenclature[C,03]{$p_{n,t,k}$}{Power consumption/production of prosumer $n$ in hour $t$ for signature $k$ [kWh]}
\nomenclature[C,04]{$p^{\mathrm{im}}_t$}{Energy imported into the community in hour $t$ [kWh]}
\nomenclature[C,05]{$p^{\mathrm{ex}}_t$}{Energy exported from the community in hour $t$ [kWh]}
\nomenclature[C,06]{$p^{\mathrm{pen}}_t$}{Energy exceeding the capacity limit agreed upon at the connection point between the community and the distribution grid in hour $t$ [kWh]}
\nomenclature[C,07]{$l_{n,t,k}$}{Residential base-load consumption for prosumer $n$ in hour $t$ for signature $k$ [kWh]}
\nomenclature[C,08]{$b_{n,t,k}$}{Battery charging (positive) or discharging (negative) for prosumer $n$ in hour $t$ for signature $k$ [kW]}
\nomenclature[C,09]{$e_{n,t,k}$}{Battery state of energy for prosumer $n$ in hour $t$ for signature $k$ [kWh]}
\nomenclature[C,10]{$p^{\mathrm{TCL}}_{n,t,k}$}{Power usage of heat pump for prosumer $n$ in hour $t$ for signature $k$ [kW]}
\nomenclature[C,11]{$\tau_{n,t,k}$}{Interior temperature for prosumer $n$ in hour $t$ for signature $k$ [$^\circ$C]}
\nomenclature[C,12]{$ev_{n,t,k}$}{Electric vehicle charging (positive) or discharging (negative) for prosumer $n$ in hour $t$ for signature $k$ [kW]}
\nomenclature[C,13]{$s_{n,t,k}$}{Electric vehicle state of energy for prosumer $n$ in hour $t$ for signature $k$ [kWh]}
\nomenclature[C,14]{$\lambda^{(.)}$}{Dual variable for all lower-level equality constraints [-]. Indexed as appropriate.}
\nomenclature[C,15]{$\mu^{(.)}$}{Dual variable for all lower-level inequality constraints [-]. Indexed as appropriate.}
\printnomenclature

\end{document}